\documentclass[review, sort&compress, numbers]{elsarticle}
\usepackage[fleqn]{amsmath}
\usepackage{graphicx}
\usepackage{lineno,hyperref}
\usepackage{latexsym}
\usepackage{amssymb}
\usepackage{amsmath}
\usepackage[center]{subfigure}
\usepackage{float}
\usepackage{color}
\usepackage{geometry}
\usepackage{diagbox}
\usepackage[section]{placeins}
\modulolinenumbers[5]

\graphicspath{{./figs/},{./logo/}}

\journal{}

\begin{document}
\begin{frontmatter}
\title{Deep neural network methods for solving forward and inverse problems of time fractional diffusion equations with conformable derivative}

\author[mymainaddress]{Yinlin~Ye}
\author[mymainaddress]{Yajing~Li\corref{mycorrespondingauthor}}
\ead{hliyajing@163.com}
\author[mymainaddress]{Hongtao~Fan}
\author[mymainaddress]{Xinyi~Liu}
\author[mymainaddress]{Hongbing~Zhang}

\cortext[mycorrespondingauthor]{Corresponding author}
\address[mymainaddress]{College of Science, Northwest A\&F University, Yangling, Shaanxi 712100, China}

\begin{abstract}
Physics-informed neural networks (PINNs) show great advantages in solving partial differential equations. In this paper, we for the first time propose to study conformable time fractional diffusion equations by using PINNs. By solving the supervise learning task, we design a new spatio-temporal function approximator with high data efficiency. L-BFGS algorithm is used to optimize our loss function, and back propagation algorithm is used to update our parameters to give our numerical solutions. For the forward problem, we can take IC/BCs as the data, and use PINN to solve the corresponding partial differential equation. Three numerical examples are are carried out to demonstrate the effectiveness of our methods. In particular, when the order of the conformable fractional derivative $\alpha$ tends to $1$, a class of weighted PINNs is introduced to overcome the accuracy degradation caused by the singularity of solutions. For the inverse problem, we use the data obtained to train the neural network, and the estimation of parameter $\lambda$ in the equation is elaborated. Similarly, we give three numerical examples to show that our method can accurately identify the parameters, even if the training data is corrupted with 1\% uncorrelated noise. 
\end{abstract}


\end{frontmatter}


\section{Introduction}

Fractional diffusion equations are often used to simulate complex phenomena in many fields, such as mechanics of wake-up or creep in polymer systems \cite{metzler00}, kinematics in viscoelastic media \cite{mainardi10}, solute transport in porous media with fractured geometry \cite{benson00}, etc. The superior modeling ability of fractional diffusion equation has aroused great interest in numerically solving such problems. Many excellent numerical methods are produced. Such methods include finite difference method\cite{zhuxu19}, finite element method\cite{karaa18}, spectral method\cite{samiee19} and so on. Although these numerical algorithms have high accuracy, most of them are time-consuming in grid generation, and can not perfectly integrate the actual data into the existing algorithms. 
\par In the past few years, machine learning algorithms have developed rapidly, especially neural network methods. The results of using neural networks to solve partial differential equations have sprung up. For example, Deep Ritz method \cite{e18} uses an energy minimization formulation as the loss function of artificial neural networks to solve Partial differential equations. Raissi et al. \cite{maziar19} introduced a physical information neural network (PINN) that is trained to solve supervised learning tasks about any given physical laws described by general nonlinear partial differential equations. When solving the forward and inverse problems of partial differential equations, they put forward two methods of continuous form and discrete form, which obtained good solving accuracy, and the speed of this method was the same when solving the forward and inverse problems of partial differential equations. Subsequently, in \cite{jincai21}, they adopted PINNs to directly encode the governing equations into the deep neural network via automatic differentiation, so as to overcome some limitations for simulating incompressible laminar and turbulent flows. Sheng and Yang in \cite{shengyang21} proposed a penalty-free neural network (PFNN) method, which can effectively solve a class of second-order boundary-value problems on complex geometries. Numerical experiments showed that PFNN was superior to several existing methods in accuracy, flexibility and robustness. Recently, more and more scholars began to pay attention to using neural network methods to solve fractional partial differential equations, and relevant literatures have emerged one after another. Pang et al. \cite{pang19} extended PINNs to fractional PINNs (fPINNs) to solve space-time fractional advection-diffusion equations, and systematically studyed their convergence. A novel element of the fPINNs was the hybrid approach that they introduced for constructing the residual in the loss function using both automatic differentiation for the integer-order operators and numerical discretization for the fractional operators. Then Pang et al. in \cite{pang20} extended PINNs to parameter and function inference for integral equations such as nonlocal
Poisson and nonlocal turbulence models, and refered to them as nonlocal PINNs (nPINNs). Qu et al. \cite{qu21} proposed a neural network method based on Legendre polynomials to
solve space and time fractional diffusion equations.
\par Over the past few decades, various definitions of fractional derivative including Grunwald-Letikov fractional derivatives, Riemann-Liouville fractional derivatives, Caputo fractional derivatives, Riesz fractional derivatives, etc, have been reported in many literatures. These fractional derivatives do not satisfy the chain rule. In 2014, Khalil et al.\cite{khalil14} introduced a new fractional derivative, which performed well and followed the Leibniz rule and the chain rule, called the conformable derivative.
Due to its effectiveness and applicability, conformable derivatives have been applied to Newtonian mechanics \cite{chung15}, quantum mechanics \cite{anderson15}, arbitrary time scale problems \cite{benkhettou16}, diffusion transport \cite{iyiola17}, neutron dynamics\cite{anaya21} and other fields. Because the conformable derivatives have good properties, the research on the theory and algorithm of partial differential equations with conformable derivatives has also attracted extensive interest. \cite{cenesiz17} studied the stochastic solution of equations with conformable time derivative where the space operators may correspond to fractional Brownian motion, or a L\'{e}vy process, or a general semigroup in a Banach space, or a process killed upon exiting a bounded domain in $\mathbb{R}^d$. \cite{hyder21} introduced a modern approach for solving the nonlinear evolution equations in the frame of a recent generalized conformable derivative. A conformable delay perturbation of matrix exponential function was offered to give the representation of solutions for linear nonhomogeneous conformable fractional delay differential equations, and the existence and uniqueness of solutions and Ulam-Hyers stability of the equations were proved in \cite{mahmudov21}. The delayed exponential matrix function in conformable version was established and used to derive the expression of the solution for homogeneous and nonhomogeneous equations respectively in \cite{xiaowang21}. \cite{singh20} introduced novel approximate numerical approach, so called "extended reduced conformable differential transform method (ERC-DTM)" which was the implementation of DTM for conformable time fractional partial differential equations with proportional delay. In \cite{bayrak21}, the truncated solution of space-time fractional differential equations, including conformable derivative was constructed by the help of residual power series method. However, there are few literatures on using neural network methods to numerically solve the conformable time fractional diffusion equations.
\par 
In this paper, we consider using PINNs to solve the conformable time fractional diffusion equation in the following form,
\begin{equation}
\left\{\begin{array}{l}
T^{\alpha} u(t, x)-\lambda u_{x x}(t, x) = 0,~~ 0<\alpha<1,~ x \in \Omega, ~t \in (0,T],  \\
u(0, x)=g(x),  \\
u(t, 0)=\varphi(t),
\end{array}\right.
\label{eq1}
\end{equation}
where $u(t,x)$ represents the solution of the equation, $\lambda\in\mathbb{R}$ is a parameter in the equation, $\Omega$ is a subset of $\mathbb{R}$, $T^{\alpha}$ is the conformable derivative.
Our aim is to give the numerical solution of the conformable time fractional diffusion equation based on the parameter $\lambda$ at a given time.
\par We summarize the main contributions and findings as follows:
\par(1) The common fractional derivatives, such as Riemann Liouville derivative and Caputo derivative, do not meet the chain rule and can not directly encode the deep neural network through the chain rule, so that the corresponding fractional differential equations can not be solved by PINNs. However, we found a fractional derivative with good properties, that is, the conformable derivative, which satisfies Leibniz's rule and chain rule. Therefore, in this paper, the PINN method is used to solve the fractional diffusion equation with integrated derivative for the first time, and experiments show that our method has good simulation effect. According to this, we have filled in the defect that neural network can not be used to solve fractional differential equations in some previous articles directly.
\par(2) As $\alpha\rightarrow 1$, due to the influence of the conformable derivative, the singularity of the solution of the equation \eqref{eq1} increases, resulting in the decline of the accuracy of our method. Therefore, we propose a weighted PINN method to deal with the influence of the singular solution by constraining the solution of the equation \eqref{eq1}.
\par(3) We use the PINN method to solve the inverse problem of the equation \eqref{eq1} and accurately identify the parameters of the equation \eqref{eq1}. Even if the training data is corrupted with 1\% uncorrelated Gaussian noise, our prediction is still robust. 
\par The rest of this paper is arranged as follows. Section 2 briefly introduces the definition and related properties of the conformable derivative.
In section 3, we present the main idea of our neural network method for solving the conformable time-fractional diffusion equation. In section 4, we solve the forward problem of the equation, i.e., taking IC/BCs as data to train our predictive solutions by minimizing the loss function. We solved the inverse problem, that is, using the data we obtained, to predict the parameter $\lambda$ of the equation by using the neural network method in section 5. Finally, we summarize the work of this paper and look forward to the future work in section 6.

\section{Preliminaries}

\par This section gives some basic definitions and properties regarding the conformable derivative. 

\noindent{\bf{Definition 1.}}
	Given a function $f$: $[0,\infty) \to \mathbb{R}$ and $t \textgreater 0$, for all $t \textgreater 0$, $\alpha \in (0,1)$, then the conformable fractional derivative of $f$ of order $\alpha$ is defined by
	\begin{equation}
	T_{\alpha}(f)(t)=\lim _{\varepsilon \rightarrow 0} \frac{f(t+\varepsilon t^{1-\alpha})-f(t)}{\varepsilon}. 
	\end{equation}
	If $f$ is $\alpha$-differentiable in some $(0,a)$, $a \textgreater 0$, and $\lim _{t \rightarrow 0+}f^{(\alpha)}(t)$ exists, then define,
	$$f^{(\alpha)}(0)=\lim _{t \rightarrow 0+}f^{(\alpha)}(t). $$
	In addition, if the conformable fractional derivative of $f$ of order $\alpha$ exists, then we simply say $f$ is $\alpha-$differentiable.
	\par Noted that $T_ {\alpha}(t^{p}) = p t^{p - \alpha}$is established. Further, the definition of conformable derivative is consistent with the classical Riemann-Liouville fractional derivative and the Caputo fractional derivative of polynomials, up to a constant multiple.

\noindent{\bf{Definition 2.}}
	If a function $f$: $[0,\infty) \to \mathbb{R}$ is $\alpha-$differentiable at $t_{0} \textgreater 0$, $\alpha \in (0,1]$, then $f$ is continuous at $t_{0}$.

\noindent{\bf{Definition 3.}}
	If $\alpha \in (0,1]$ and $f, g$ is $\alpha-$differentiable at $t \textgreater 0$, then, \\
	(1) $T_{\alpha}(a f+b g)=a T_{\alpha}(f)+b T_{\alpha}(g)$, for all $a, b \in \mathbb{R}. $\\
	(2) $T_{\alpha}\left(t^{p}\right)=p t^{p-\alpha}$ for all $p \in \mathbb{R}. $\\
	(3) $T_{\alpha}(\lambda)=0$, for all constant functions $f(t)=\lambda. $\\
	(4) $T_{\alpha}(f g)=f T_{\alpha}(g)+g T_{\alpha}(f). $\\
	(5) $T_{\alpha}\left(\frac{f}{g}\right)=\frac{g T_{\alpha}(f)-f T_{\alpha}(g)}{g^{2}}. $\\
	(6) In addition, if $f$ is differentiable, then $T_{\alpha}(f)(t)=t^{1-\alpha} \frac{d f}{d t}(t). $

\section{Methodology}

For the conformable time-fractional diffusion equations\eqref{eq1}, we use $f(t,x)$ to determine the left part, i.e.,
\begin{equation}
f(t,x) = T^{\alpha} u(t, x)-\lambda u_{x x}(t, x),
\end{equation}
where $u(t,x)$ can be approximated by using a deep neural networks. Different from other definitions of fractional derivatives, conformable derivatives satisfy the chain rule, so we can approach it by automatic differentiation. To highlight the simplicity of implementing this idea, we use a Python snippet containing TensorFlow. TensorFlow is a symbolic mathematics system based on data flow programming, which is widely applied in the programming implementation of various machine learning algorithms. Thus, $u(t,x)$ can be simply defined as:
$$
\begin{aligned}
\textbf{def} \quad & u(t , x): \\
& u \  = \  neural \_ net ( tf . concat ([ t , x] ,1) , \  weights , \  biases ) \\
& \textbf{return} \quad u \\
\end{aligned}
$$
Similarly, $f(t, x)$ can be simply defined as:
$$
\begin{aligned}
\textbf{def} \quad & f (t , x):\\
& u \  = \  u(t , x)\\
& u_t \  = \  tf . gradients (u, t )[0]\\
& u_x \  = \  tf . gradients (u, x)[0]\\
& u_{xx} \  = \  tf . gradients (u _ x , x)[0]\\
& f \  = \  u_t \  - \  \lambda*u_{xx}\\
& \textbf{return} \quad f\\
\end{aligned}
$$
\par During the training of PINN, we approximate the solution of the equation by minimizing the following loss function. 
\begin{equation}
Loss = MSE_{u}+MSE_{f},
\label{eq4}
\end{equation}
where MSE represents the mean square error and the parameters of neural networks $u(t,x)$ and $f(t,x)$ are shared. The loss function consists of two parts, which are specifically expressed as follows: 
$$
MSE_{u} = \frac{1}{N_{u}} \sum_{i=1}^{N_{u}}\left|u\left(t_{u}^{i}, x_{u}^{i}\right) - u^{i}\right|^{2}, 
$$
where $\{t_{u}^{i},x_{u}^{i},u^{i}\}_{i=1}^{N_{u}}, N_{u}$ represent the initial and the boundary training data of the equation, the number of training points for ICs and BCs, respectively, 
$$
MSE_{f} = \frac{1}{N_{f}} \sum_{i=1}^{N_{f}}\left|f\left(t_{f}^{i}, x_{f}^{i}\right)\right|^{2}
$$
where $\{{t_{f}^{i}, x_{f}^{i}}\}$ represent the collocation of neural network $f(t,x)$. The error $MSE_{u}$ corresponds to  the initial values and boundary values of the equation, while the error $MSE_{f}$ is used to reinforce the structure imposed by equation \eqref{eq1} within a finite set of collocation points. 
\par The neural network $u(t,x)$ trained by IC/BCs and the neural network $f(t,x)$ containing physical information both contribute to the loss function and share the hyper-parameters, which together constitute the PINN architecture for conformable time-fractional diffusion equations as displayed in the following figure \ref{pinns}.
\begin{figure}[htb]
	\centering
	\includegraphics[scale=0.5]{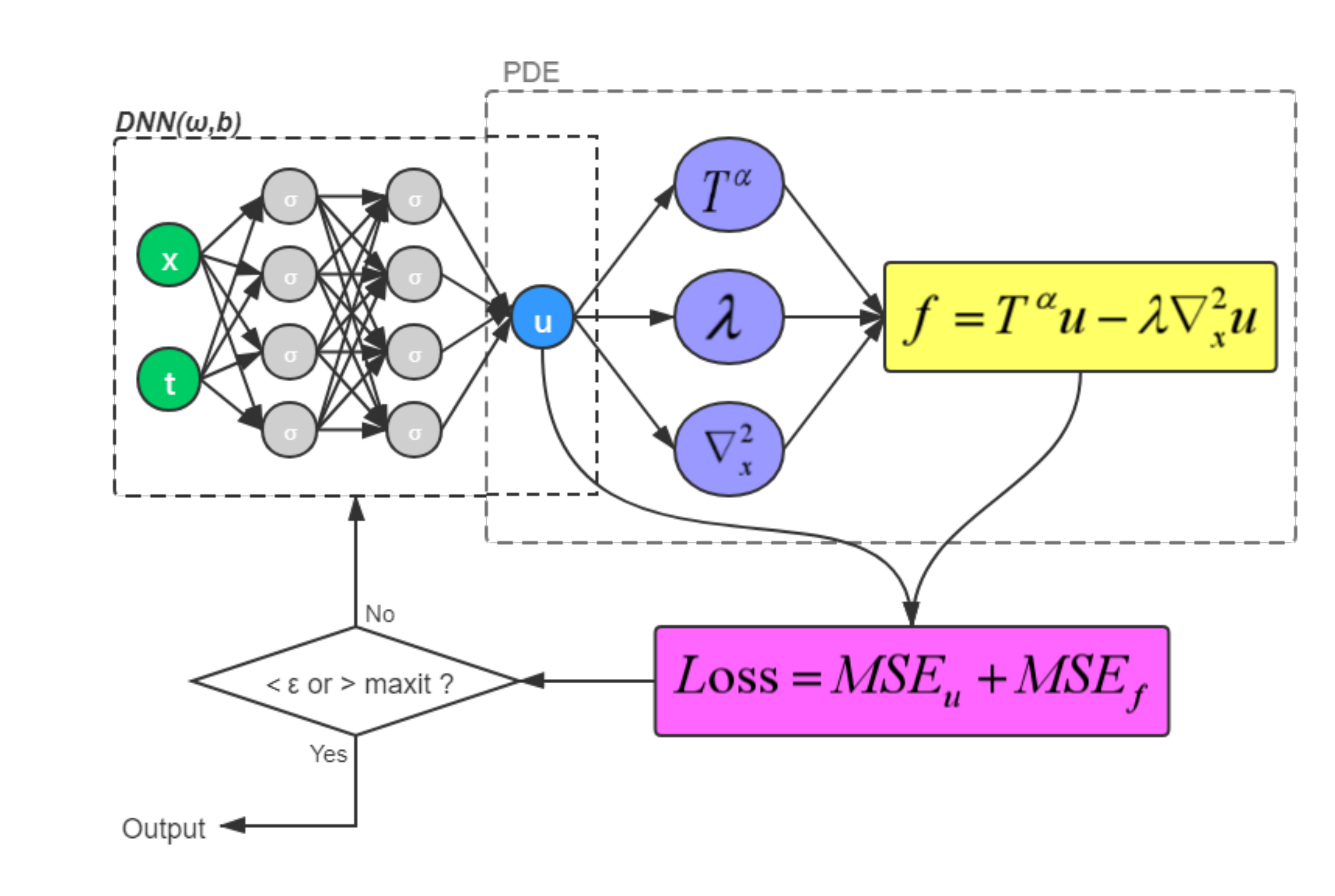}
	\caption{Schematic of PINN for the conformable time-fractional diffusion equations. }
	\label{pinns}
\end{figure}
\par By the employment of PINN, we can take IC/BCs as the data to solve the corresponding partial differential equation since it is obvious that the conformable time-fractional diffusion equation is wellposed, with suitable IC/BCs. This is referred to as the forward problem, and we shall report our results in section 4; On the other hand, assuming that we have some of the training data, which may contain unrrelated noise, the parameterized equations are learned by making use of it. Likewise, it is referred to as the inverse problem and in section 5.

\section{Forward problems}

In this section, we aim to solve the forward problem for the conformable time-fractional diffusion equations, that is, to obtain the neural network approximation of the solution of the equation using IC/BCs and prior physical information. It is worth mentioning that for wellposed forward problems, fractional partial differential equations can usually be uniquely solved by using IC/BCs. 
\par For this forward problem, by taking IC/BCs as data the loss function \eqref{eq4} can be established and designed as
\begin{equation}
	Loss = MSE_{u} + MSE_{f} = MSE_{IC} + MSE_{BC} + MSE_{f}
	\label{eq5}
\end{equation}
where
$$
MSE_{IC} = \frac{1}{N_{IC}} \sum_{i=1}^{N_{IC}}\left|u\left(0, x_{IC}^{i}\right) - g(x_{IC}^{i})\right|^{2}, 
$$
and
$$
MSE_{BC} = \frac{1}{N_{BC}} \sum_{i=j}^{N_{BC}}\left|u\left(t_{BC}^{j}, 0\right) - \phi(t_{BC}^{j}) \right|^{2},
$$
here $\{x_{IC}^{i},g(x_{IC}^{i})\}_{i=1}^{N_{IC}}, \{t_{BC}^{j},\phi(t_{BC}^{j})\}_{j=1}^{N_{BC}}, N_{IC}, N_{BC}$ denote the initial value of the equation, the boundary training data of the equation, the number of training points for ICs and the number of training points for BCs, respectively. Besides, $MSE_{f}$ applied to represent prior physical information can be rewritten as
$$
\begin{aligned}
MSE_{f} & = \frac{1}{N_{f}} \sum_{i=1}^{N_{f}}\left|f\left(t_{f}^{i}, x_{f}^{i}\right)\right|^{2} \\
& = \frac{1}{N_{f}} \sum_{i=1}^{N_{f}}\left| T^{\alpha} u(t_{f}^{i}, x_{f}^{i})-\lambda u_{x x}(t_{f}^{i}, x_{f}^{i}) \right|^{2}
\end{aligned}
$$
where $\{{t_{f}^{i}, x_{f}^{i}}\}$ indicate the collocation of neural network $f(t,x)$ and $N_{f}$ represents the number of collocation points. The error $MSE_{f}$ in the neural network $f(t,x)$ is applied to reinforce the structure imposed by equation \eqref{eq1} within a finite set of collocation points. 
\par Next, we will carry out three examples to demonstrate the applicability of our neural network method in solving conformable time-fractional diffusion equations. All experiments are implemented on a computer, which is configured as follows: Intel(R) Core(TM) i5-8265U CPU @ 1.60 GHz 1.80 GHz, using Python3.6 + Tensorflow1.15. The same configuration is applied in the inverse problem in Section 5.

\begin{figure}[htb]
	\centering
	\includegraphics[scale=0.4]{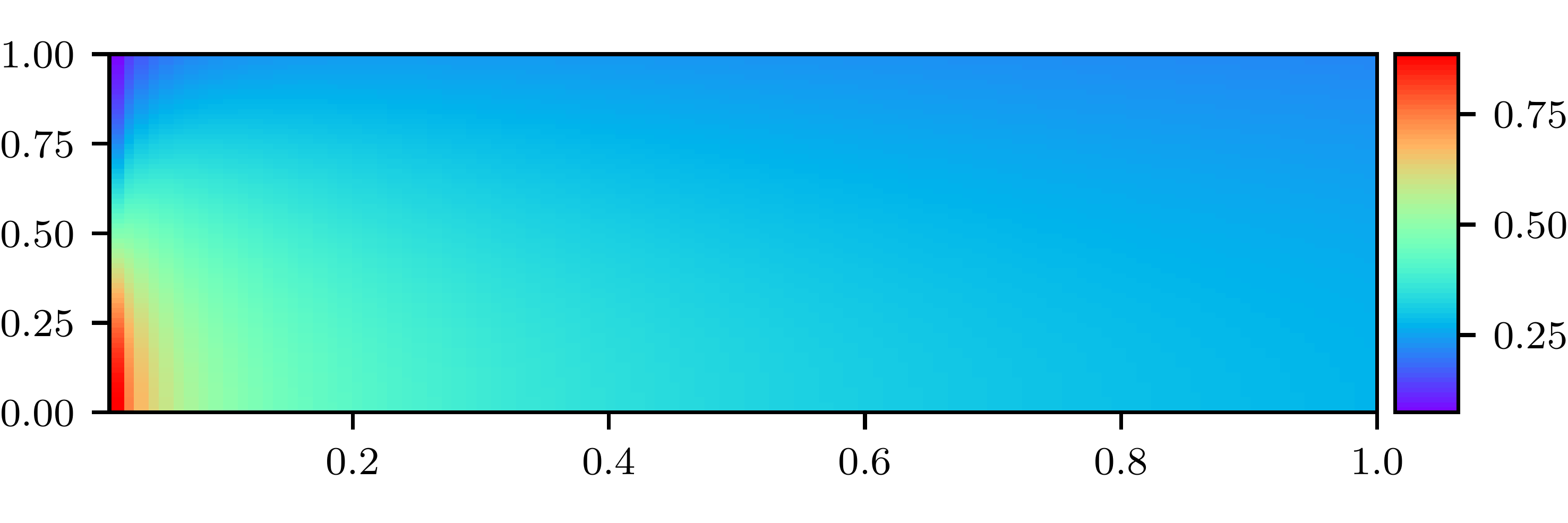}
	\caption{Analytical solution of the equation when $\alpha = 0.5$.}
	\label{fig1}
\end{figure}

\noindent{\bf{Example 1.}}
We consider equation \eqref{eq1} with $\alpha = 0.5 $, $\lambda = 0.5073 $, which has an analytical solution of the following form,
$$
u(t,x)=\sqrt{\frac{\alpha}{4 \pi \lambda t^{\alpha}}} \exp \left\{-\frac{\alpha}{4 \lambda t^{\alpha}} x^{2}\right\},
$$
the image of the analytical solution is exhibited in figure \ref{fig1}.
\par
\begin{figure}[htb]
	\centering
	\includegraphics[scale=1]{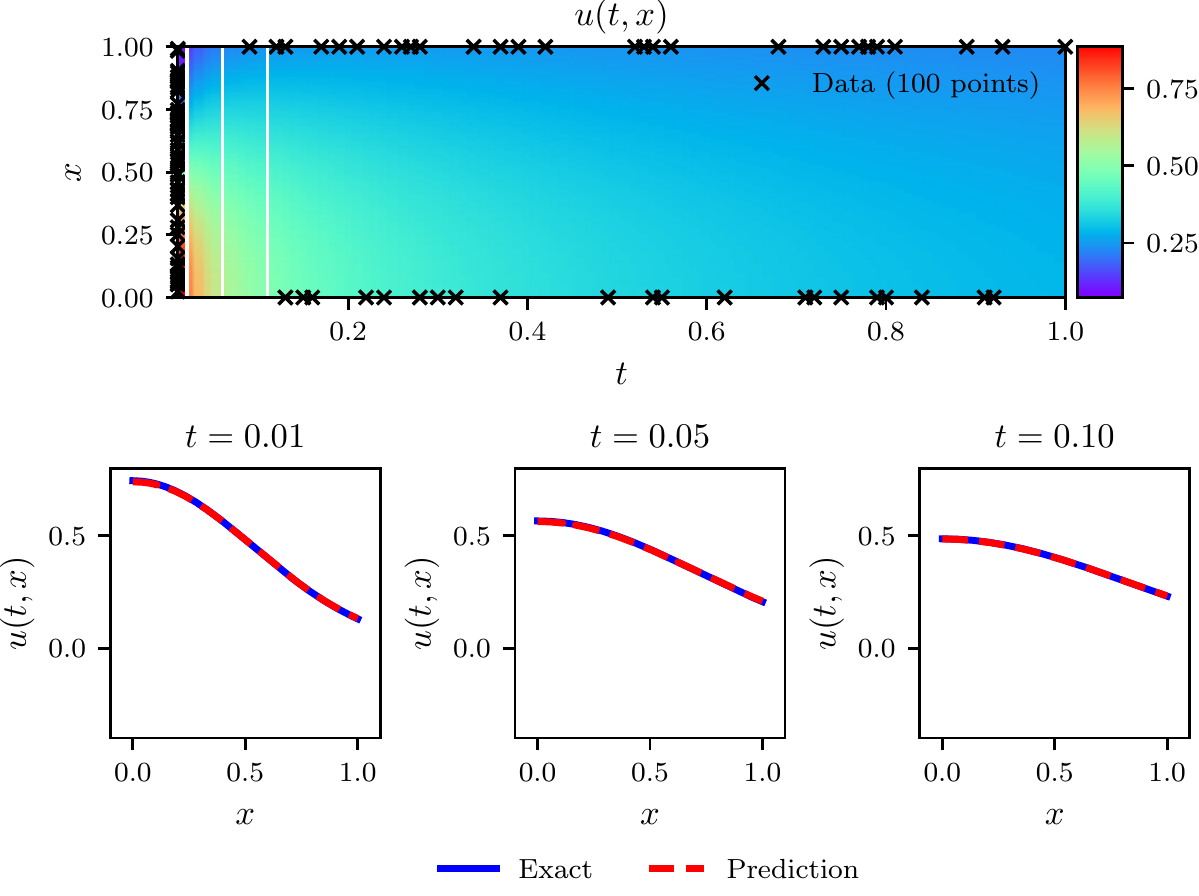}
	\caption{Comparison between the exact solution of the equation at $\alpha = 0.5$ and the exact solution and the predicted solution at three different times of $t = 0.01, 0.05, 0.10$.}
	\label{fig2}
\end{figure}

\begin{figure}[htb]
	\centering
	\includegraphics[scale=1]{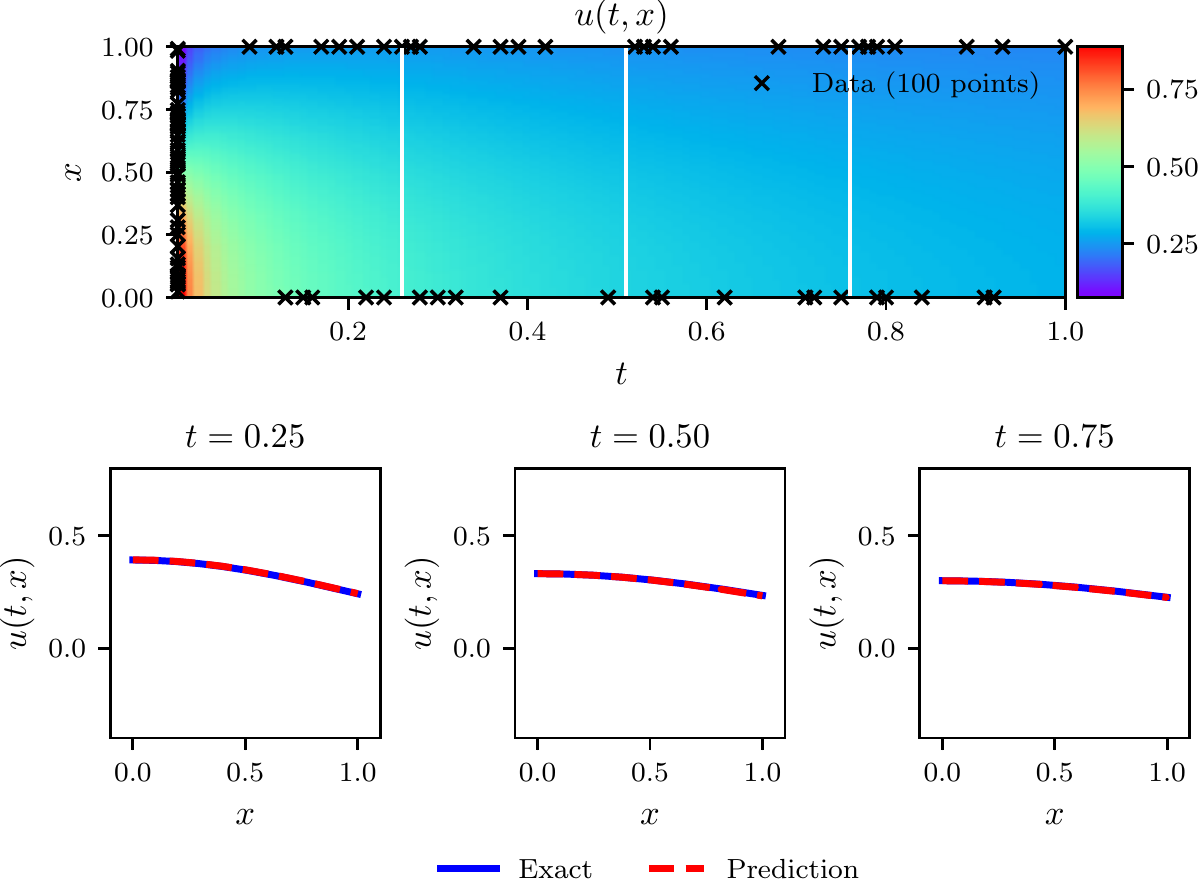}
	\caption{Comparison between the exact solution of the equation at $\alpha = 0.5$ and the exact solution and the predicted solution at three different times of $t = 0.25, 0.50, 0.75$.}
	\label{fig3}
\end{figure}

 Figures \ref{fig2}-\ref{fig3} reveal the predicted results of our neural network method for data-driven solutions of the conformable time-fractional diffusion equation. Specifically, given a set of initial and boundary data with $N_{u} = N_{IC} + N_{BC} = 100$ and a set of collocation points with $N_{f} = 10000$, both of them are randomly distributed. Next, we employ the mean square error loss function defined in equation \eqref{eq4} to train 3021 parameters of the deep neural network that has 9 hidden layers with 20 neurons at each layer, and learn the solution $u(t,x)$ of the conformable time-fractional diffusion equation. The hyperbolic tangent function serve as our activation function over here. At the top pannel in figures \ref{fig2}-\ref{fig3}, we display the space-time solution $u(t,x)$ predicted by our neural network and the location of the initial and boundary training data. What here must stress specially is that, all classical numerical methods for solving partial differential equations require discretization of the equations, but our method does not require any discretization in time or space. In the case, the exact solution given is analytically available and relative $\mathbb{L}_{2}$ error of the predicted solution is measured as 2 $\cdot 10^{-3}$. At the bottom panel of of figures \ref{fig2}-\ref{fig3}, a more detailed evaluation of the predicted solution is offered. To be specific, figure \ref{fig2} is presented with the comparison between the exact solution and the predicted solution at three different moments of $t=0.01,0.05,0.10$. In addition, the comparison between the exact solution and the predicted solution $t=0.25,0.50,0.75$ is displayed in figure \ref{fig3}.  
\par
\begin{figure}[htb]
	\centering
	\includegraphics[scale=0.4]{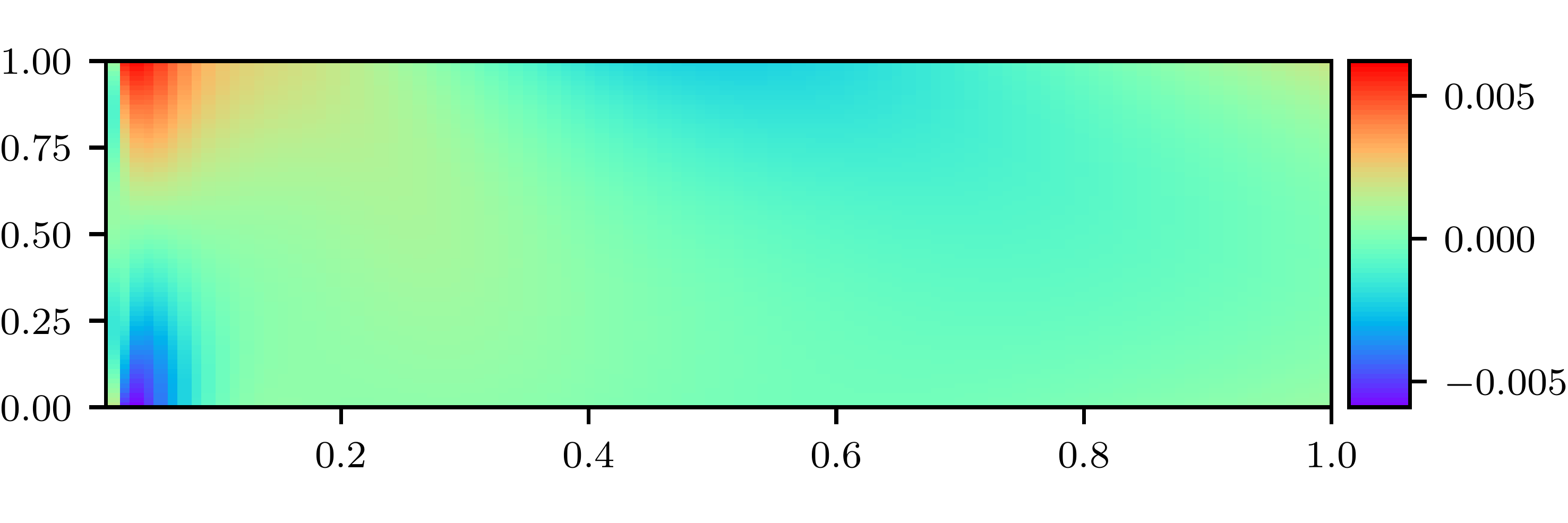}
	\caption{Error between the predicted solutions and exact solutions of the equation when $\alpha = 0.5$.}
	\label{fig4}
\end{figure}

\begin{figure}[htb]
	\centering
	\includegraphics[scale=0.4]{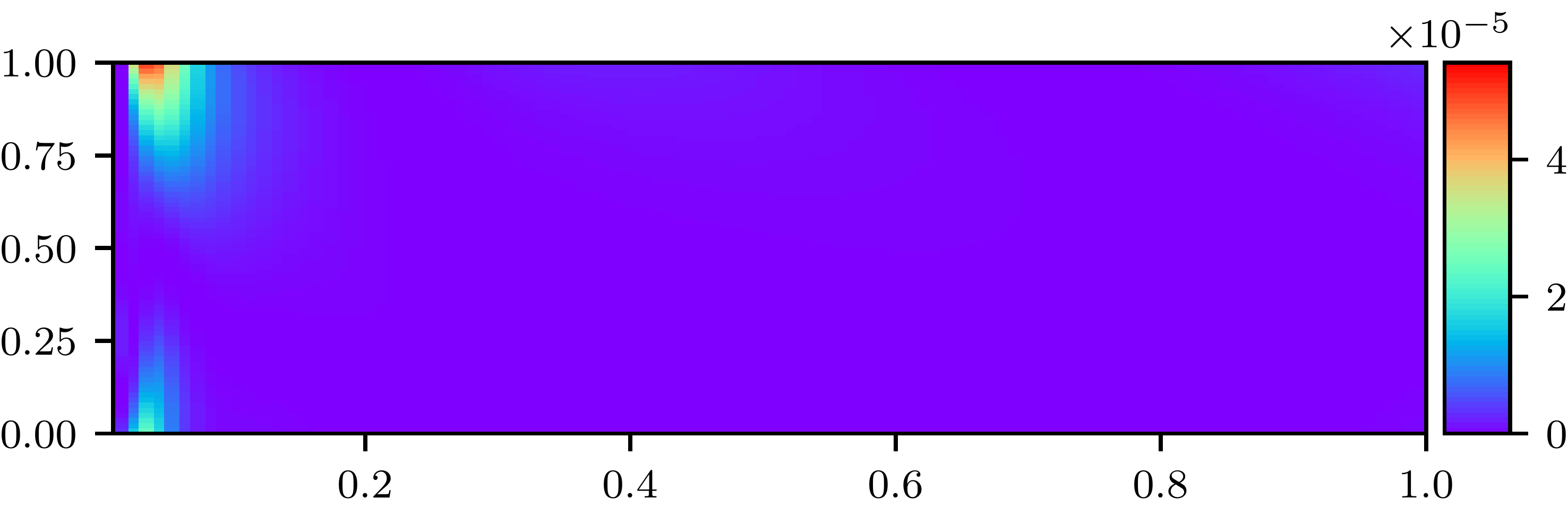}
	\caption{Variance of the predicted solutions and exact solutions of the equation when $\alpha = 0.5$.}
	\label{fig5}
\end{figure}

 The error between predicted solution and exact solutions when $\alpha = 0.5$ is drawn in figure \ref{fig4}. We can observe that the error is particularly close to 0, while the average error is about $2.0 \cdot 10^{-3}$. In figure \ref{fig5}, we exhibit the variances between predicted solution and exact solutions for $\alpha = 0.5$. It is worth mentioning that the mean square error is about $4.2 \cdot 10^{-7}$. From figures \ref{fig4}-\ref{fig5}, we can find that the predicted solution is a good approximation of the exact solutions for the conformable time-fractional diffusion equation at $\alpha = 0.5$. 

\begin{figure}[htb]
	\centering
	\includegraphics[scale=0.4]{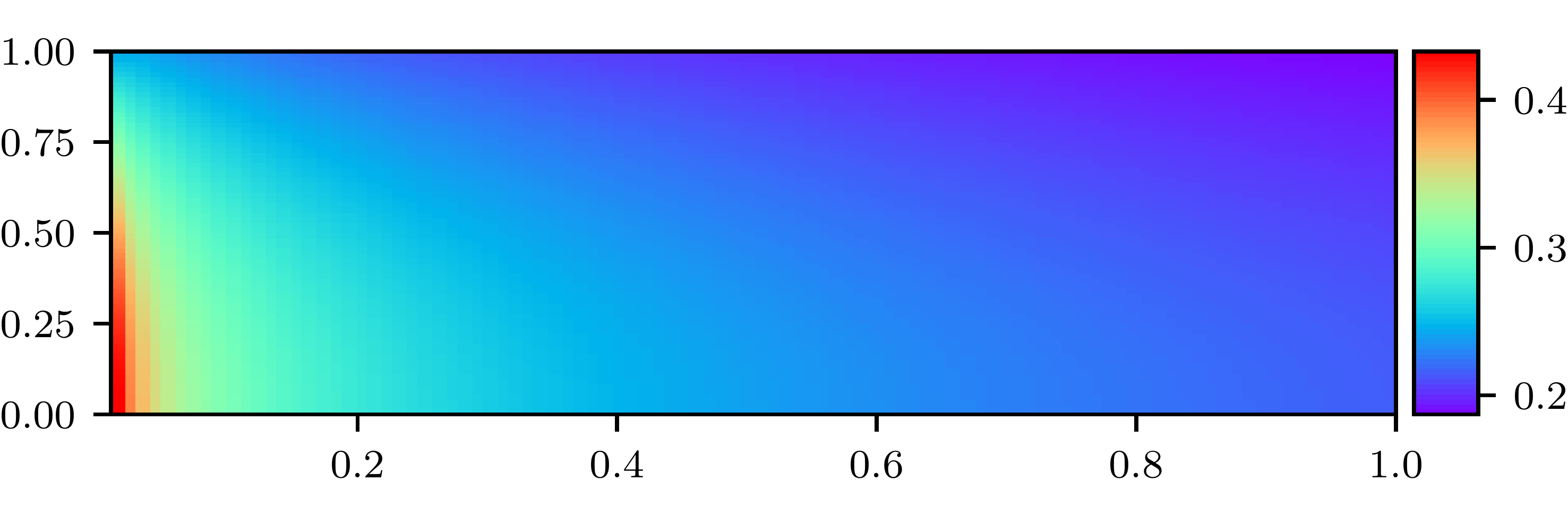}
	\caption{Analytical solution of the equation when $\alpha = 0.3$.}
	\label{fig6}
\end{figure}

\noindent{\bf{Example 2.}}
In this example, the analytical solution of equation \eqref{eq1} with $\alpha = 0.3$, $\lambda = 0.5073$ is given by
$$
u(t,x)=\sqrt{\frac{\alpha}{4 \pi \lambda t^{\alpha}}} \exp \left\{-\frac{\alpha}{4 \lambda t^{\alpha}} x^{2}\right\},
$$
whose image is described in \ref{fig6}. 
\par
\begin{figure}[htb]
	\centering
	\includegraphics[scale=1]{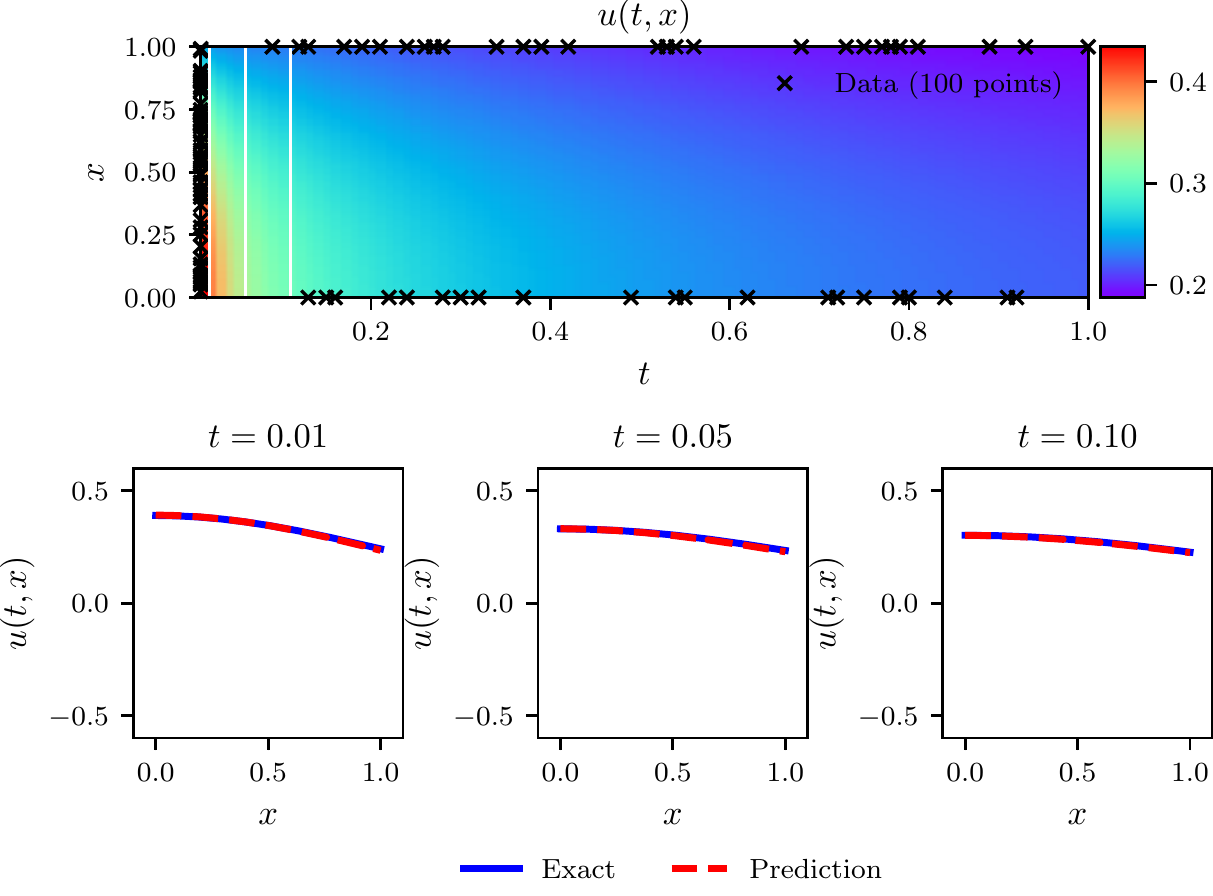}
	\caption{Comparison between the exact solution of the equation at $\alpha = 0.3$ and the exact solution and the predicted solution at three different times of $t = 0.01, 0.05, 0.10$.}
	\label{fig7}
\end{figure}

\begin{figure}[htb]
	\centering
	\includegraphics[scale=1]{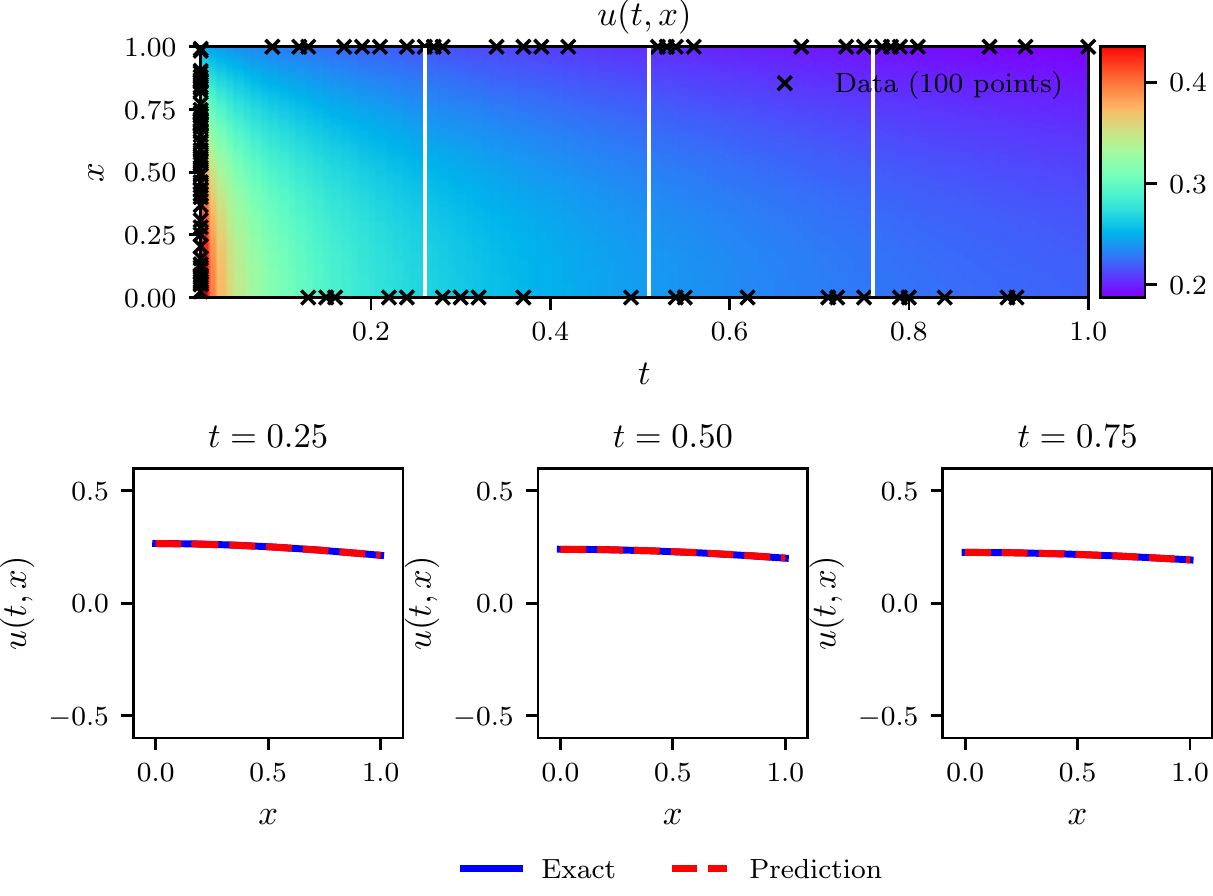}
	\caption{Comparison between the exact solution of the equation at $\alpha = 0.3$ and the exact solution and the predicted solution at three different times of $t = 0.25, 0.50, 0.75$.}
	\label{fig8}
\end{figure}

 For the conformable time-fractional diffusion equation, the results predicted by neural network methodin are listed in figures \ref{fig7}-\ref{fig8}. Concretely, both a set of initial and boundary data with $N_{u} = N_{IC} + N_{BC} = 100$ and a set of collocation points with $N_{f} = 10000$ are randomly distributed. Then we make use of the mean square error loss function defined in equation \eqref{eq4} to learn the solution $u(t,x)$ of the conformable time-fractional diffusion equation by training 3021 parameters while the deep neural network has 9 hidden layers with 20 neurons per layer. Here, we pick the hyperbolic tangent function as activation function. Figures \ref{fig7}-\ref{fig8} manifest the space-time solution $u(t,x)$ predicted by our neural network and the location of the initial and boundary training data at the top pannel. For this problem, the relative error of the predicted solution is measured as 1.8 $\cdot10^{-3}$ in the $\mathbb{L}_{2}$ norm. At the bottom panel of of figures \ref{fig7}-\ref{fig8}, we reveal the comparison between the exact solution and the predicted solution at $t=0.01,0.05,0.10$ and $t=0.25,0.50,0.75$, respectively. 
\par
\begin{figure}[htb]
	\centering
	\includegraphics[scale=0.4]{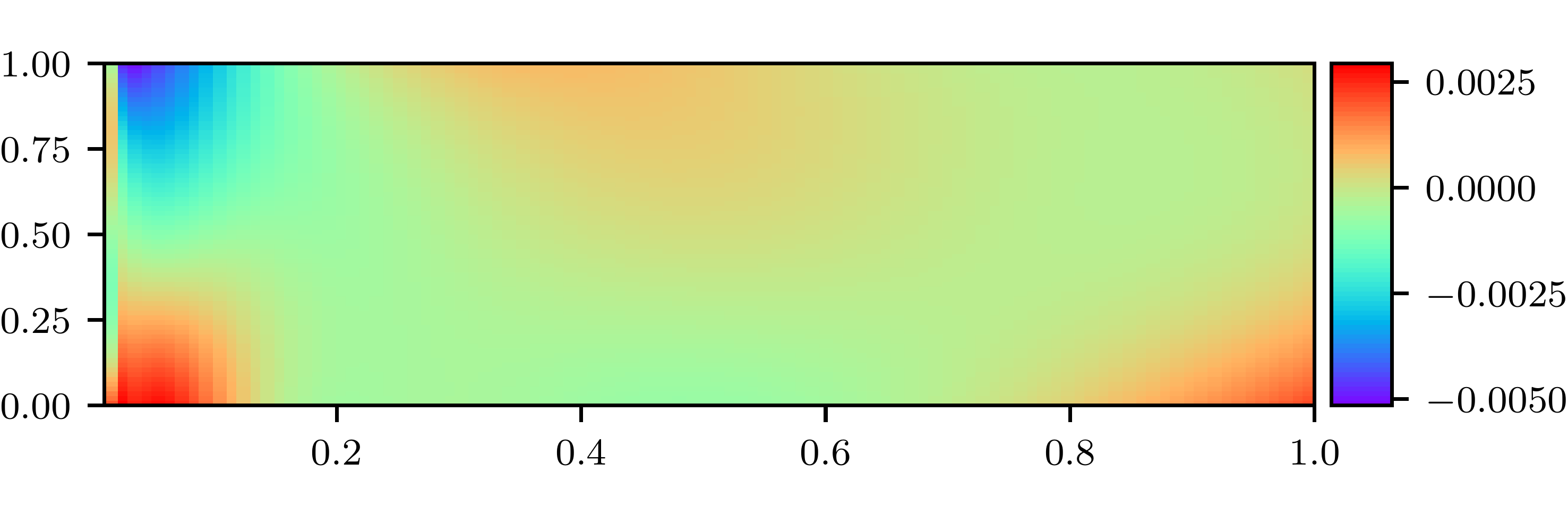}
	\caption{Error between the predicted solutions and exact solutions of the equation when $\alpha = 0.3$.}
	\label{fig9}
\end{figure}

\begin{figure}[htb]
	\centering
	\includegraphics[scale=0.4]{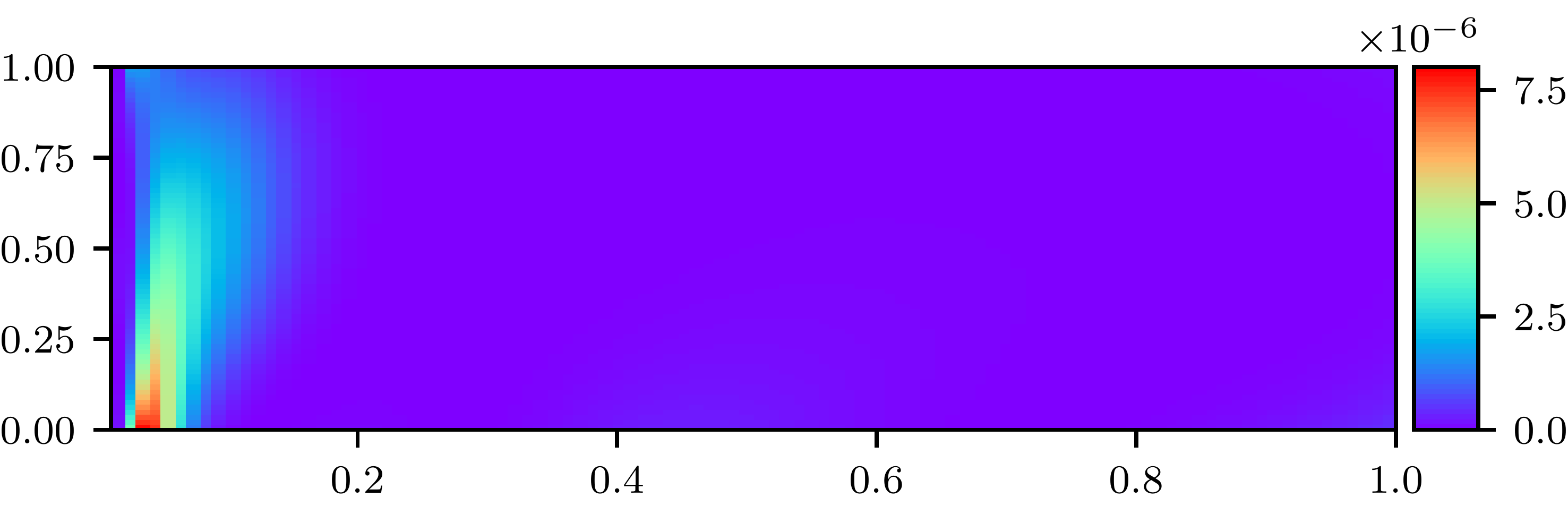}
	\caption{Variance of the predicted solutions and exact solutions of the equation when $\alpha = 0.3$.}
	\label{fig10}
\end{figure}

 In figure \ref{fig9}, we exhibit the error of predicted solution and exact solutions with $\alpha = 0.3$. What we can see is that the error is close to 0 particularly and the average error is about $1.8 \cdot 10^{-3}$. Besides, the variances between predicted solution and exact solutions with $\alpha = 0.3$ are portrayed in figure \ref{fig10}, while the mean square error is about $2.6 \cdot 10^{-7}$. In addition, for the conformable time-fractional diffusion equation at $\alpha = 0.3$, we can discover that the predicted solution approximate the exact solutions pretty well from figures \ref{fig9}-\ref{fig10}. 

\begin{figure}[htb]
	\centering
	\includegraphics[scale=0.4]{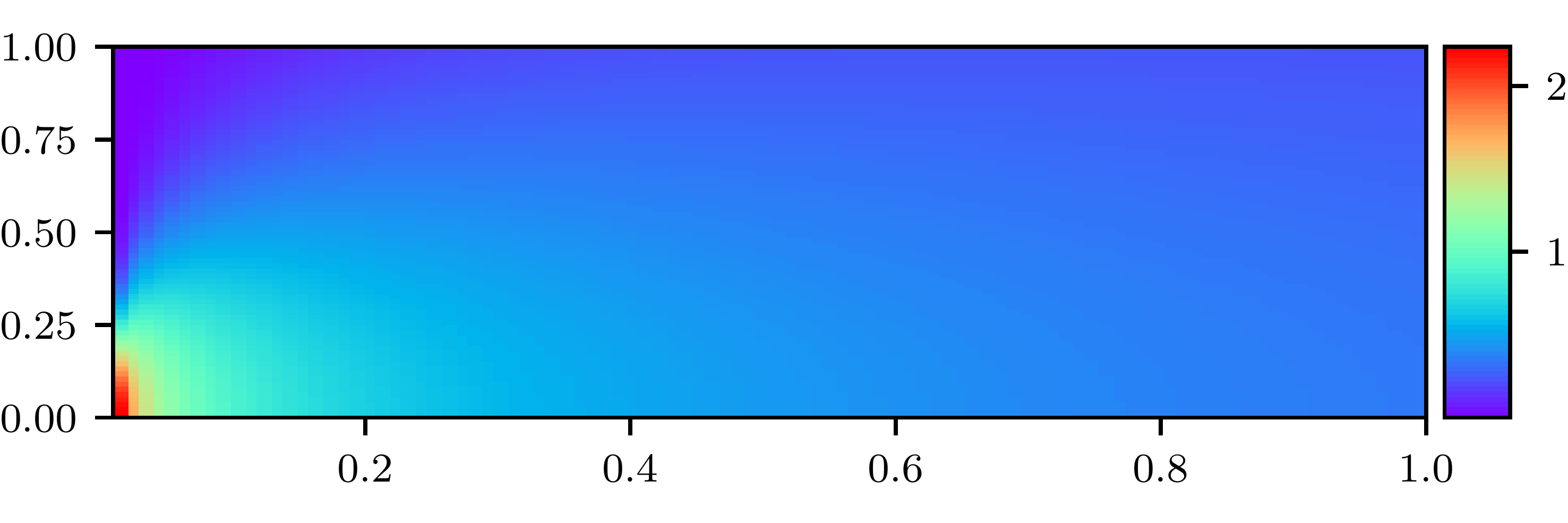}
	\caption{Analytical solution of the equation when $\alpha = 0.8$.}
	\label{fig11}
\end{figure}

\noindent{\bf{Example 3.}}
Figure \ref{fig11} reveals the analytical solution of equation \eqref{eq1}, which can be represented as the following form,
$$
u(t,x)=\sqrt{\frac{\alpha}{4 \pi \lambda t^{\alpha}}} \exp \left\{-\frac{\alpha}{4 \lambda t^{\alpha}} x^{2}\right\},
$$
where $\alpha = 0.8$, $\lambda = 0.5073$.
\par
\begin{figure}[htb]
	\centering
	\includegraphics[scale=1]{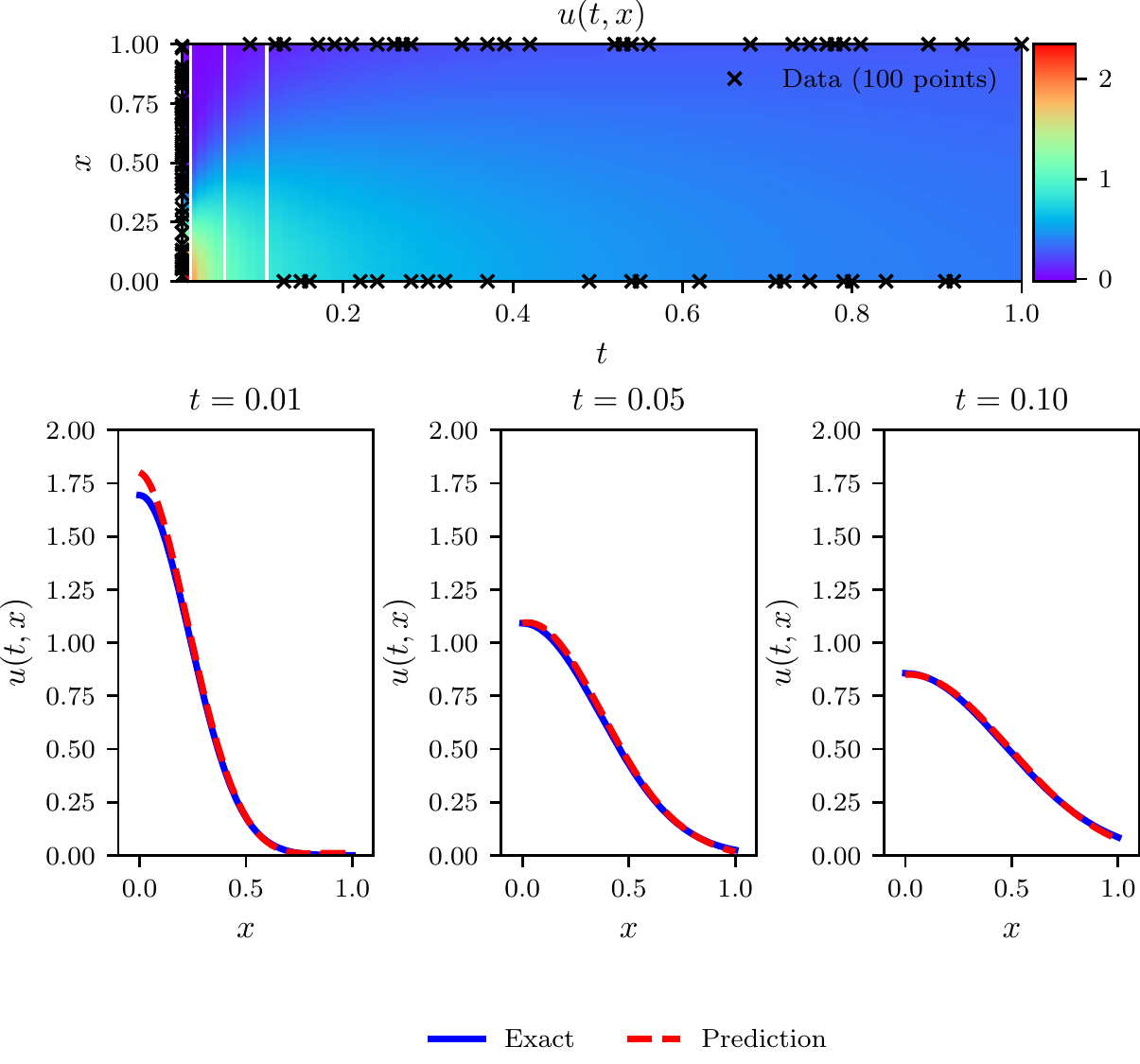}
	\caption{Comparison between the exact solution of the equation at $\alpha = 0.8$ and the exact solution and the predicted solution at three different times of $t = 0.01, 0.05, 0.10$.}
	\label{fig12}
\end{figure}

\begin{figure}[htb]
	\centering
	\includegraphics[scale=1]{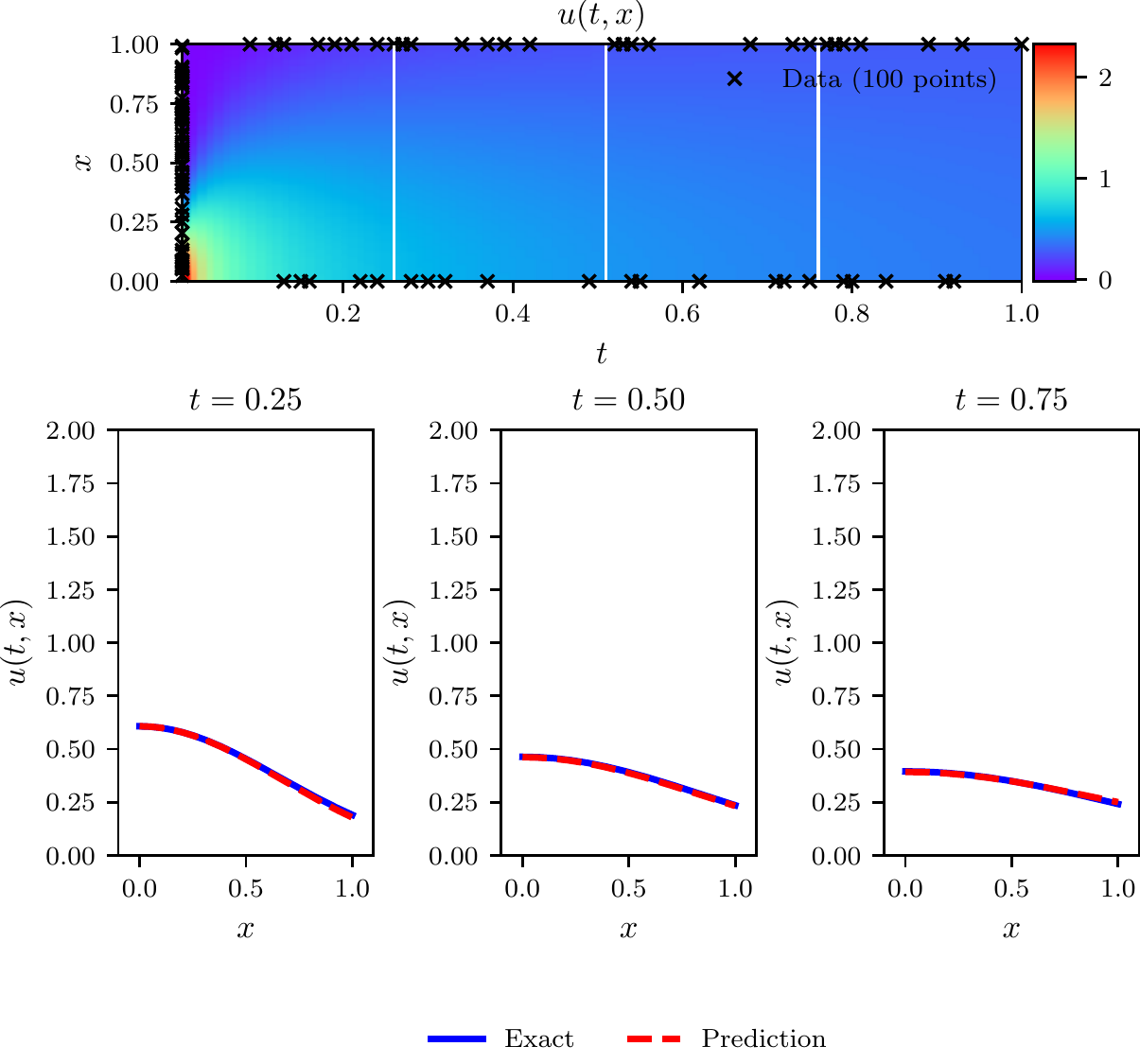}
	\caption{Comparison between the exact solution of the equation at $\alpha = 0.8$ and the exact solution and the predicted solution at three different times of $t = 0.25, 0.50, 0.75$.}
	\label{fig13}
\end{figure}

 Figures \ref{fig12}-\ref{fig13} provide the predicted results by making use of our neural network method for the conformable time-fractional diffusion equation. Specifically, both the initial and boundary data with $N_{u} = N_{IC} + N_{BC} = 100$ and the collocation points with $N_{f} = 10000$ are selected randomly. And then the mean square error loss function defined in equation \eqref{eq4} is employed to train 3021 parameters of the deep neural network which has 9 hidden layers with 20 neurons at each layer and learn the solution $u(t,x)$ of the conformable time-fractional diffusion equation. In the same way, the hyperbolic tangent function is chosen as our activation function. At the top pannel of figures \ref{fig12}-\ref{fig13}, we exhibit the space-time solution $u(t,x)$ predicted by our neural network and the location of the initial and boundary training data. In addition, the exact solution given is truely available and the relative $\mathbb{L}_{2}$ error between predicted solution and exact solution is measured as 1.4 $\cdot10^{-2}$. In addition, more detailed evaluation of the predicted solution is given at the bottom panel of figures \ref{fig12}-\ref{fig13}. In detail, we give the comparison between the exact solution and the predicted solution at three different moments of $t=0.01,0.05,0.10$ in figure \ref{fig12}. Unfortunately, we can see that the predicted solution has some error at $t=0.01$ and some flaws at both $t=0.05$ and $t=0.10$ in this figure. We will discuss why this is the case and how to solve it later in this section. On the contrary, in figure \ref{fig13}, the predicted solution we got at three different moments of $t=0.25,0.50,0.75$ is perfoming well. 
\par
\begin{figure}[htb]
	\centering
	\includegraphics[scale=0.4]{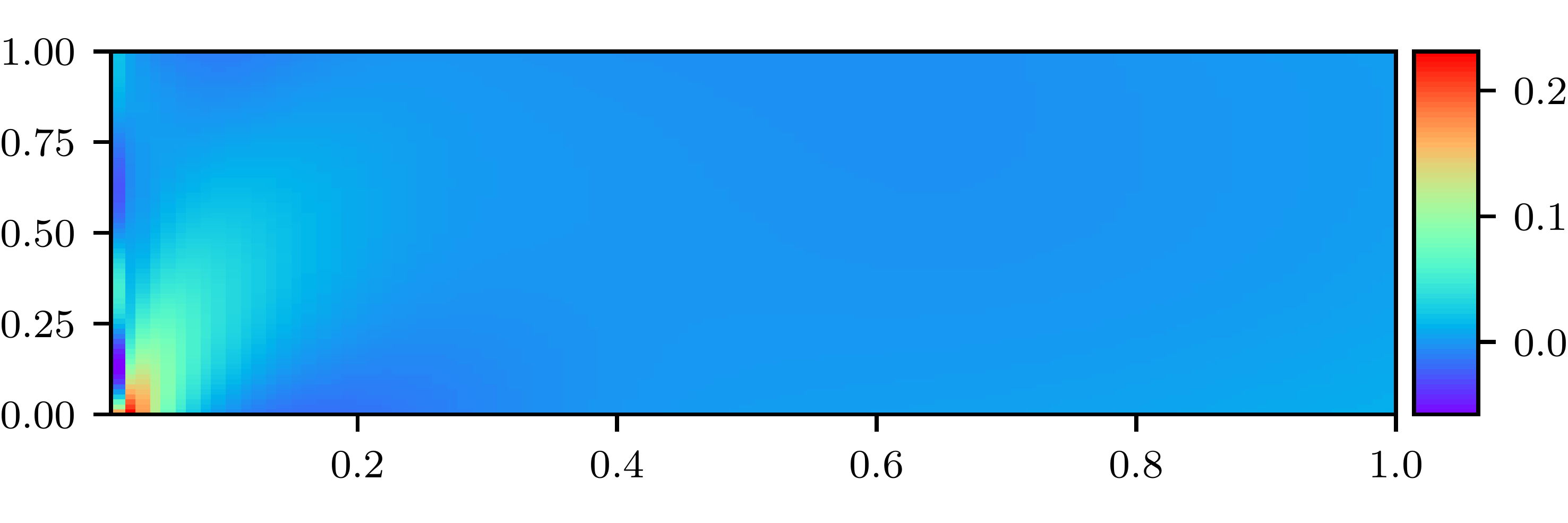}
	\caption{Error between the predicted and exact solutions of the equation when $\alpha = 0.8$.}
	\label{fig14}
\end{figure}

\begin{figure}[htb]
	\centering
	\includegraphics[scale=0.4]{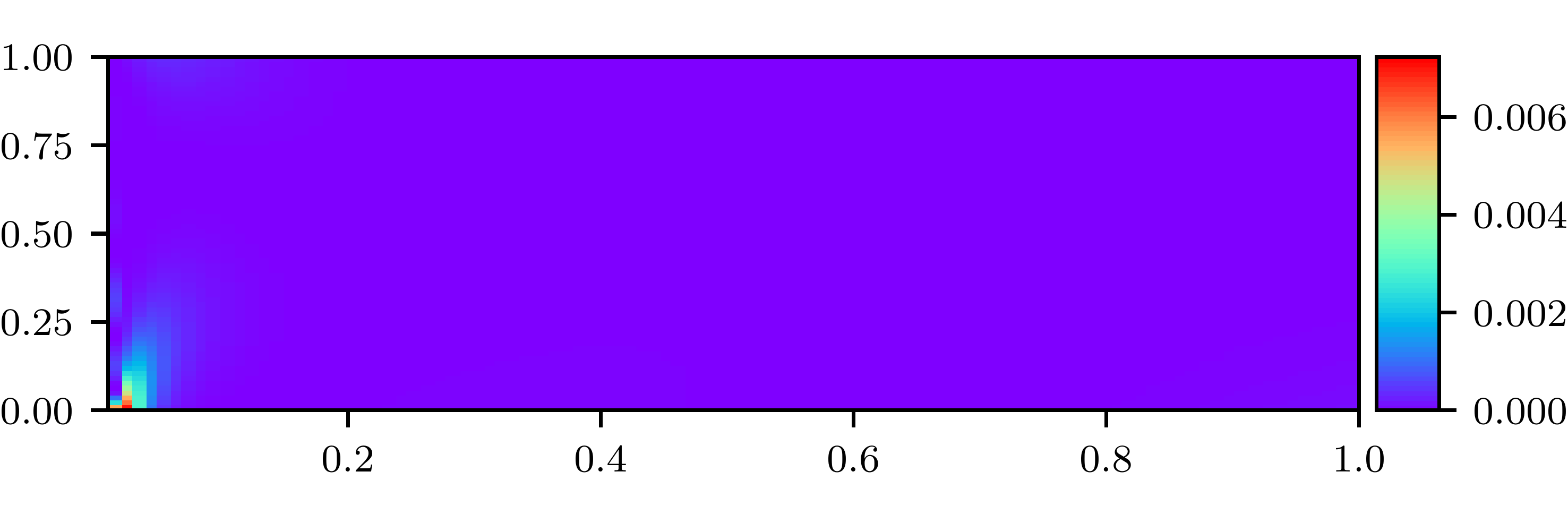}
	\caption{Variance of the predicted solutions and exact solutions of the equation when $\alpha = 0.8$.}
	\label{fig15}
\end{figure}

 In figure \ref{fig14}, the error between predicted solution and exact solutions with $\alpha = 0.8$ is described. What we show is that the error is extremly close to 0, and the average error is about $1.4 \cdot 10^{-2}$. In addition, we exhibit the variances of  predicted solution and exact solutions at $\alpha = 0.8$ in figure \ref{fig15} and the mean square error is about $3.7\cdot10^{-5}$. However, variances we got are relatively large in the singular solution region while particularly close to 0 in all other regions. For the equation at $\alpha=0.8$, the predicted solution has a good approximation out of the singular solution region in figures \ref{fig14}-\ref{fig15}.

\par Table \ref{table1} presents the mean error and mean variance between the predicted solution and the exact solution of the confoemable time-fractional diffusion equation when $\alpha=0.3,0.5,0.8$. It can be observed that, with the increase of the fractional order $\alpha$ of the conformable derivative, the singularity of the solution of the equation increases as it tends to 0. Also, the error of the solution of the equation increases as well. Therefore, we propose a weighted neural network method to deal with the complex nonlinear behavior caused by the singular solution of the equation when the fractional order $\alpha$ tends to 1.

\begin{table}
	\caption{The mean error and mean variance between the predicted solutions and exact solutions of the conformable time-fractional diffusion equation when $\alpha=0.3,0.5,0.8$.}
	\label{table1}
	\begin{tabular}{c |c c c}
		\hline
		$\alpha$ & 0.3 & 0.5 & 0.8 \\
		\hline
		Average error & $1.8 \cdot 10^{-3}$ & $2.0 \cdot 10^{-3}$ & $1.4 \cdot 10^{-2}$ \\
		\hline
		Mean square error & $4.2 \cdot 10^{-7}$ & $2.6 \cdot 10^{-7}$ & $3.7 \cdot 10^{-5}$ \\
		\hline
	\end{tabular}
\end{table}
\par In our original neural network method, the loss function is defined as the equation \eqref{eq4}. It is a difficulty that when $\alpha$ close to 1, the complex nonlinear behavior caused by the singular solution of equation is difficult to be simulated. Therefore, we consider to weight the neural network $u (t, x)$ and $f (t, x)$ and vest the mean square error $($MSE$)$ greater weight. By doing so, the predicted solution of the neural networkmethod can be constrained and the influence of singular solution of the equation can be reduced better. 
\par The parameters of neural networks $u(t,x)$ and $f(t,x)$ can be learned by minimizing the mean square error $($MSE$)$loss function, which is defined as follows:
\begin{equation}
MSE = w_{u} * MSE_{u}+ w_{f} * MSE_{f},
\end{equation}
where $w_{u}$ and $w_{f}$ imply neural network $u(t,x)$ and $f(t,x)$, respectively. Since we have used this symbol $w$ to represent the weight of the neural network, here we use $w$ as the weight symbol instead to avoid confusion. Next, we provide experimental results to verify the accuracy of our weighted neural network method.
\par
\begin{figure}[htb]
	\centering
	\includegraphics[scale=1]{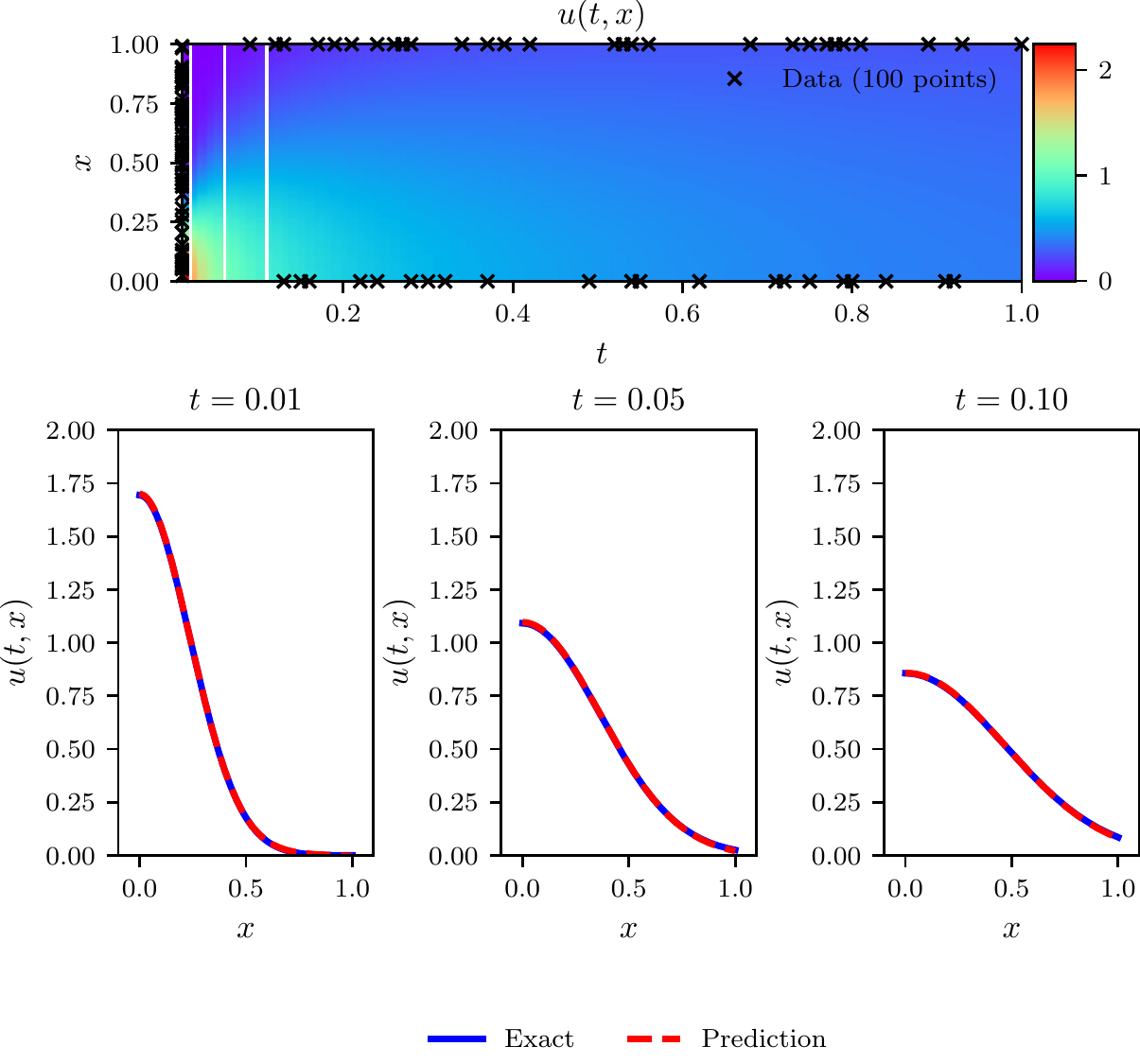}
	\caption{Comparison between the exact solution of the equation at $\alpha = 0.8$ and the exact solution and the predicted solution at three different times of $t = 0.01, 0.05, 0.10$.}
	\label{fig16}
\end{figure}

 For the equation \eqref{eq1} with $\alpha = 0.8 $, $\lambda = 0.5073 $, we take $w_{u} = 1.0, w_{f}=0.1$ in our weighted neural network method. Figure \ref{fig16} shows the predicted results of our neural network method for the data-driven solutions of the conformable time-fractional diffusion equation. To make it more concrete, given a set of initial and boundary data with $N_{u} = N_{IC} + N_{BC} = 100$ and a set of collocation points with $N_{f} = 10000$, both of them are randomly distributed. After that, we exploit the mean square error loss function defined in equation \eqref{eq4} to train 3021 parameters of the deep neural network that has 9 hidden layers with 20 neurons per layer, and learn the solution $u(t,x)$ of the conformable time-fractional diffusion equation. Similarly, the hyperbolic tangent function is selected as our activation function. 
\par Results of the prediction error is measured as 1.4 $\cdot10^{-3}$ in the relative $\mathbb{L}_{2}$ norm. The comparison between the exact solution and the predicted solution at three different times of $t=0.01,0.05$ and $0.10$ is shown, see figure \ref{fig16}. To our delight, we find that the accuracy of our weighted neural network method is fine at $t=0.01,0.05$ and $0.10$.

\par
\begin{figure}[htb]
	\centering
	\includegraphics[scale=0.4]{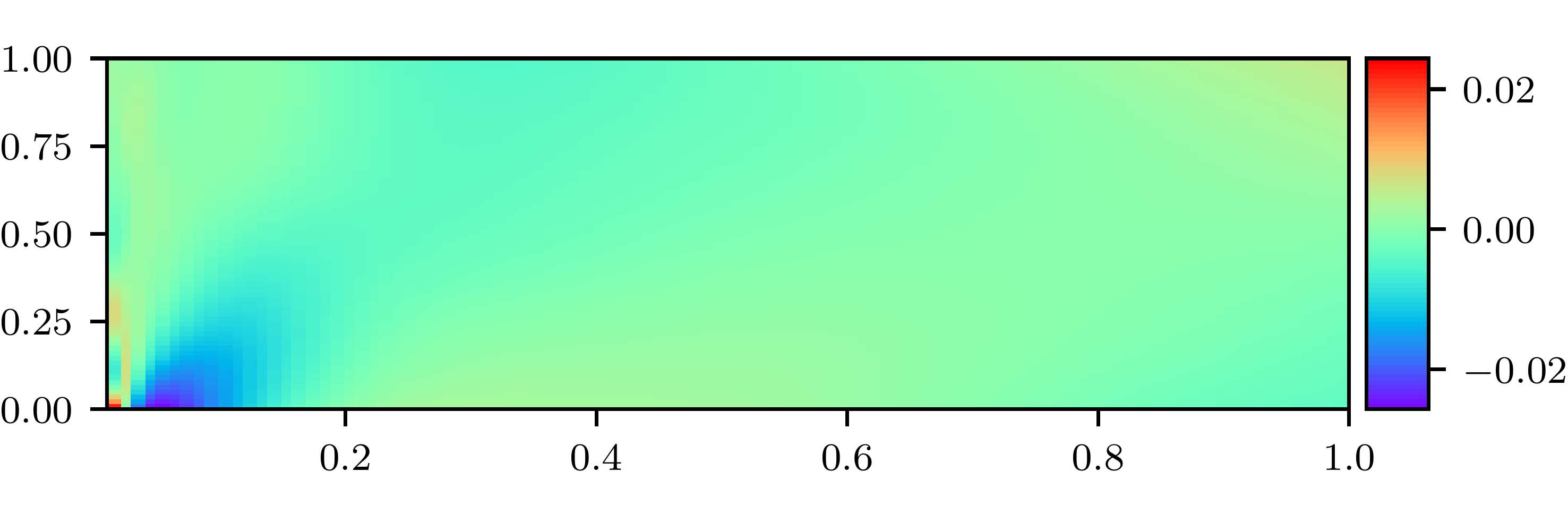}
	\caption{error of the weighted neural network predicted solutions and exact solutions of the equation when $\alpha=0.8$.}
	\label{fig17}
\end{figure}

\begin{figure}[htb]
	\centering
	\includegraphics[scale=0.4]{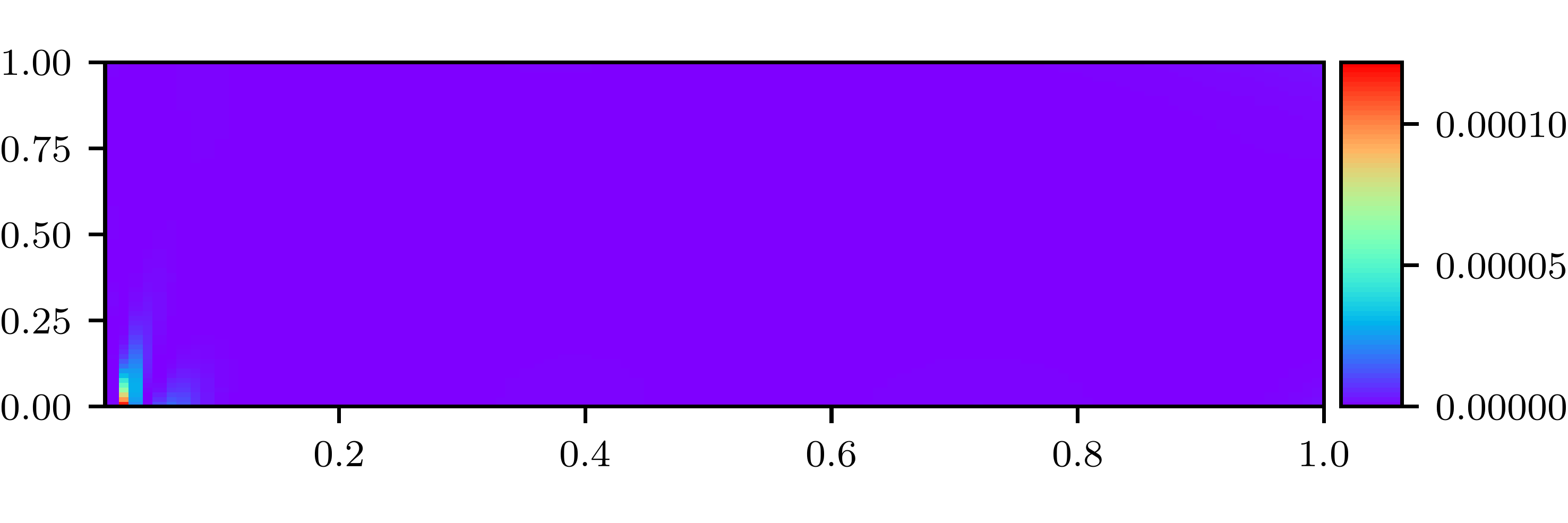}
	\caption{The variance of the weighted neural network predicted solutions and exact solutions of the equation when $\alpha=0.8$.}
	\label{fig18}
\end{figure}

 In figure \ref{fig17}, we exhibit the error between predicted solution and exact solutions when $\alpha = 0.8$. At this time, we can see that the error is adequately close to 0, and the average error is about $1.4 \cdot 10^{-3}$. Figure \ref{fig18} shows the variances of our predicted solution and exact solutions for $\alpha = 0.8$ and the mean square error is about $3.5 \cdot 10^{-7}$. Through figures \ref{fig17}-\ref{fig18}, what we want to demonstrate is that the predicted solution is a good approximation of the exact solutions for the conformable time-fractional diffusion equation at $\alpha = 0.8$. 

\par To further analyze the performance of our method, we conduct a systematic study to quantify the effects of different sampling points (training points, collocation points) and neural network structures (layers, the number of neurons in each layer) on the prediction accuracy. From Table \ref{table2}, we show the relative $\mathbb{L}_{2}-$ error generated by different number of initial and boundary training data $N_{u}$ and different number of collocation points $N_{f}$ for conformable time-fractional diffusion equation when $\alpha = 0.5$. Here we keep the structure of the 9-layer deep neural network with 20 neurons per layer unchanged. The general trend is that in the case of enough collocation points $N_{f}$, the prediction accuracy will increase with the increase of the total training data $N_{u}$. This observation highlights a key advantage of our neural network method: by configuring the point $N_{f}$ to encode the structure of the underlying physical laws, a more accurate and data-efficient learning algorithm can be obtained. 
\par Finally, Table \ref{table3} collects the relative $\mathbb {L}_{2}$ errors generated by different hidden layers and neurons in each layer for the conformable time-fractional diffusion equation when $\alpha=0.5$. Similarly, the total number of training points and collocation points is fixed at $N_{u}=100$ and $N_{f}= 10,000$, respectively. As expected, we observe that as the number of layers and the number of neurons per layer increases, this also mean that the ability of the neural network to approximate more complex functions increases, and the prediction accuracy increases accordingly.

\begin{table}[htb]
	\caption{Different initial and boundary number of training data $N_{u} $ and different collocation points  $N_{f}$, the relative $\mathbb{L}_{2}$ error between the predicted values and the exact solution $u (t, x)$. Here, the network structure is fixed at nine layers, with 20 neurons in each hidden layer.}
	\label{table2}
	\begin{tabular}{c |c c c c c}
		\hline
		\diagbox{$N_{u}$}{$N_{f}$} & 2000 & 4000 & 6000 & 8000 & 10000 \\
		\hline
		20& $1.4 \cdot 10^{-2}$ & $3.7 \cdot 10^{-2}$ & $8.5 \cdot 10^{-3}$ & $7.4 \cdot 10^{-3}$ & $5.2 \cdot 10^{-3}$ \\
		40& $1.1 \cdot 10^{-2}$ & $2.1 \cdot 10^{-2}$ & $1.1 \cdot 10^{-2}$ & $7.0 \cdot 10^{-3}$ & $5.3 \cdot 10^{-3}$ \\
		60& $9.6 \cdot 10^{-3}$ & $1.3 \cdot 10^{-2}$ & $9.6 \cdot 10^{-3}$ & $5.8 \cdot 10^{-3}$ & $4.6 \cdot 10^{-3}$ \\
		80& $8.1 \cdot 10^{-3}$ & $8.2 \cdot 10^{-3}$ & $7.2 \cdot 10^{-3}$ & $4.8 \cdot 10^{-3}$ & $4.2 \cdot 10^{-3}$ \\	
		100& $7.4 \cdot 10^{-3}$ & $6.4 \cdot 10^{-3}$ & $4.0 \cdot 10^{-3}$ & $3.6 \cdot 10^{-3}$ & $2.0 \cdot 10^{-3}$ \\
		200& $9.5 \cdot 10^{-3}$ & $4.3 \cdot 10^{-3}$ & $3.8 \cdot 10^{-3}$ & $2.6 \cdot 10^{-3}$ & $1.9 \cdot 10^{-3}$ \\
		\hline
	\end{tabular}
\end{table}

\section{Inverse problems}
In this section, we turn our attention to the inverse problem for the conformable time-fractional diffusion equations. Assume in advance that we can get some training data, at this time we want to use it to learn parameterized equations. Usually, the data we derived is mixed with some uncorrelated noise, and we want to know whether our method is still applicable in this case. 
\begin{table}[htb]
	\caption{Relative $\mathbb{L}_{2}$ error between the predicted result and the exact solution $u(t, x)$ under different hidden layers and neurons per layer. Here, the total number of training points and matching points is fixed at $N_{u}$ = 100 and $N_{f}$ = 10,000, respectively.}
	\label{table3}
	\begin{tabular}{c |c c c}
		\hline
		\diagbox{Layers}{Neurons} & 10 & 20 & 40  \\
		\hline
		2& $7.4 \cdot 10^{-2}$ & $6.8 \cdot 10^{-2}$ & $9.8 \cdot 10^{-3}$ \\
		4& $3.1 \cdot 10^{-2}$ & $7.1 \cdot 10^{-3}$ & $6.3 \cdot 10^{-3}$ \\
		6& $4.4 \cdot 10^{-2}$ & $8.8 \cdot 10^{-3}$ & $4.5 \cdot 10^{-3}$ \\
		8& $8.2 \cdot 10^{-3}$ & $4.9 \cdot 10^{-3}$ & $2.0 \cdot 10^{-3}$ \\	
		\hline
	\end{tabular}
\end{table}
$$
MSE_{data} = \frac{1}{N_{data}} \sum_{i=1}^{N_{data}}\left|u\left(0, x_{data}^{i}\right) - g(x_{data}^{i})\right|^{2}, 
$$
and
$$
\begin{aligned}
MSE_{f} & = \frac{1}{N_{data}} \sum_{i=1}^{N_{data}}\left|f\left(t_{data}^{i}, x_{data}^{i}\right)\right|^{2} \\
& = \frac{1}{N_{data}} \sum_{i=1}^{N_{data}}\left| T^{\alpha} u(t_{data}^{i}, x_{data}^{i})-\lambda u_{x x}(t_{data}^{i}, x_{data}^{i}) \right|^{2}, \\
\end{aligned}
$$
where $\{x_{data}^{i},g(x_{data}^{i})\}_{i=1}^{N_{data}}, N_{data}$ signify the training data on $u(t,x)$ and the number of training points, respectively. In particular, the loss $MSE_{u}$ corresponds to the training data on u(t, x) while $MSE_{f}$ enforces the structure imposed by the conformable time-fractional diffusion equation at a finite set of collocation points, whose number and location are taken to be the same as the training data. At this time, $\lambda$ is the parameter we want to learn.
\par Similarly, we bestow three examples to demonstrate the applicability of our neural network method in solving the inverse problem of conformable time-fractional diffusion equations. What is worth mentioning is that all experiments are performed on a computer, whose configuration is the same as in Section 4. 

\noindent{\bf{Example 4.}}
In this example, we consider equation \eqref{eq1} with $\alpha = 0.5$. In order to illustrate the effectiveness of our method, we create a training dataset by randomly generating $N = 2000$ points across the entire spatio-temporal domain from the exact solution corresponding to $\lambda = 0.5073$. Next, we employ $N_{data} = 2000$ training points to train 3021 parameters of the deep neural network that has 9 hidden layers with 20 neurons at each layer by using the Adam optimizer and L-BFGS optimizer to minimize the mean square error loss function defined in equation \eqref{eq4} and learn the solution $u(t,x)$ of the conformable time-fractional diffusion equation. Here, the hyperbolic tangent function served as activation function. After training, the network is calibrated to predict the entire solution u(t, x), as well as the unknown parameter $\lambda$. 
\begin{figure}[htb]
	\centering
	\includegraphics[scale=1]{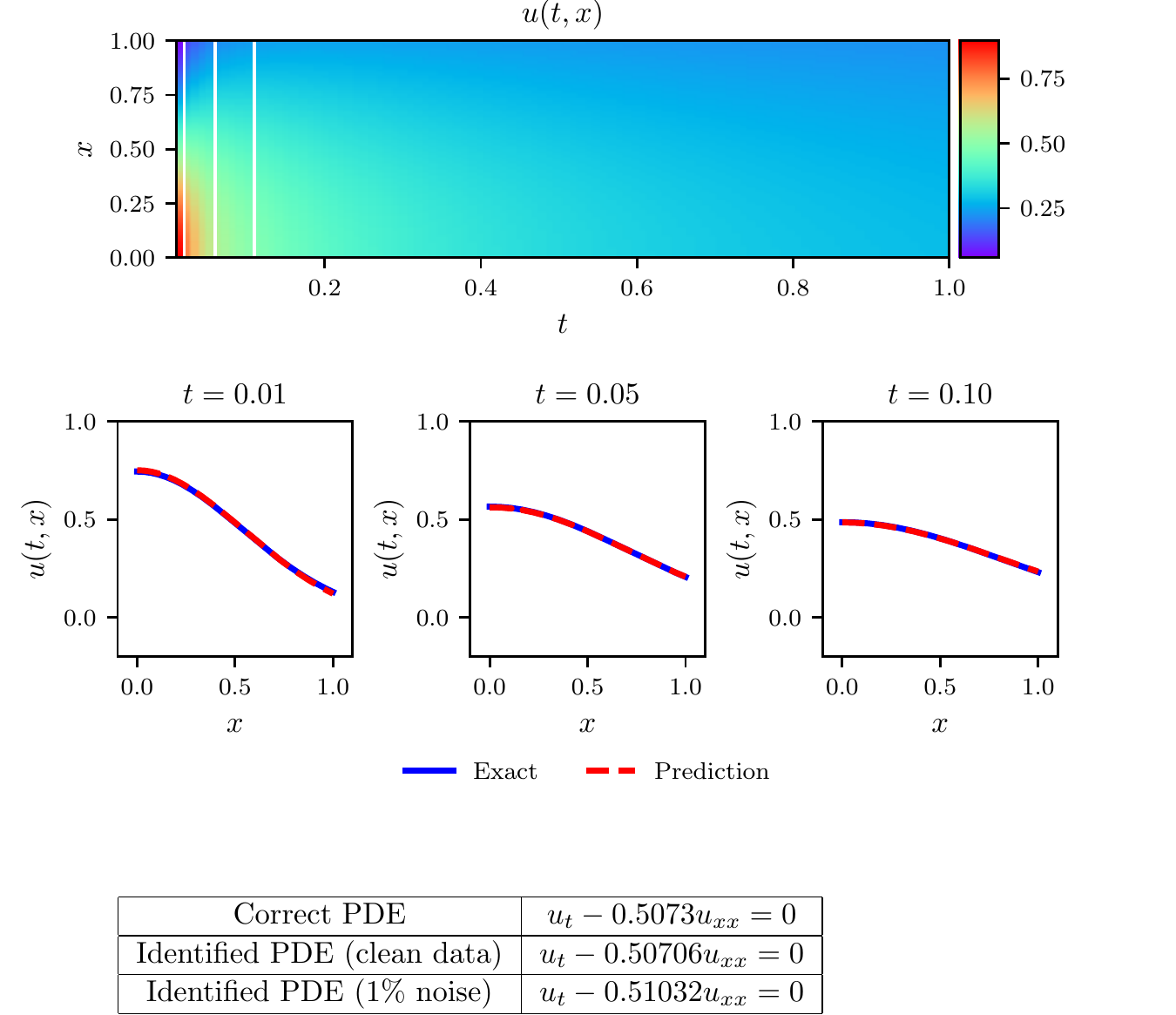}
	\caption{Equation with $\alpha = 0.5$: Top: Solution u(t, x) along with the temporal locations of the three training snapshots. Middle: Training data and exact solution corresponding to the three temporal snapshots depicted by the dashed vertical lines in the top panel. Bottom: Correct partial differential equation along with the identified one obtained by learning $\lambda$.}
	\label{fig20}
\end{figure}
\par In figure \ref{fig20} , we summarize the results for the conformable time-fractional diffusion equation with $\alpha = 0.5$. At the top panel, we supply the predicted solution of the equation whose analytical solution can be found in figure \ref{fig1}. At the middle pannel, the graph of the predicted solution and the exact solution at $t=0.01,0.05,0.10$ is drawn. At the bottom pannel, we furnish the correct PDE and the identified PDE with clean data and 1$\%$ noise. We observe that our neural network method can accurately identify the unknown parameter $\lambda$, even if the training data is corrupted with noise. Specifically, the estimation errors of $\lambda$ is 0.05$\%$ in the case of clean training data. Although the training data is corrupted with 1$\%$ uncorrelated Gaussian noise, the prediction is still robust, returning the error of $\lambda$ is 0.60$\%$.

\noindent{\bf{Example 5.}}
The equation \eqref{eq1} with $\alpha = 0.3$ is taken into account in this example. In order to illustrate the effectiveness of our method, we create a training dataset by generating $N = 2000$ points  randomly across the entire spatio-temporal domain from the exact solution corresponding to $\lambda = 0.5073$. And then, we make use of $N_{data} = 2000$ training points to train 3021 parameters of the deep neural network which has 9 hidden layers with 20 neurons per layer by using the Adam optimizer and L-BFGS optimizer to minimize the mean square error loss function defined in equation \eqref{eq4} and learn the solution $u(t,x)$ of the conformable time-fractional diffusion equation. The hyperbolic tangent function is singled out as our activation function over here. Later, the network is calibrated to predict the entire solution u(t, x) and the unknown parameter $\lambda$. 
\begin{figure}[htb]
	\centering
	\includegraphics[scale=1]{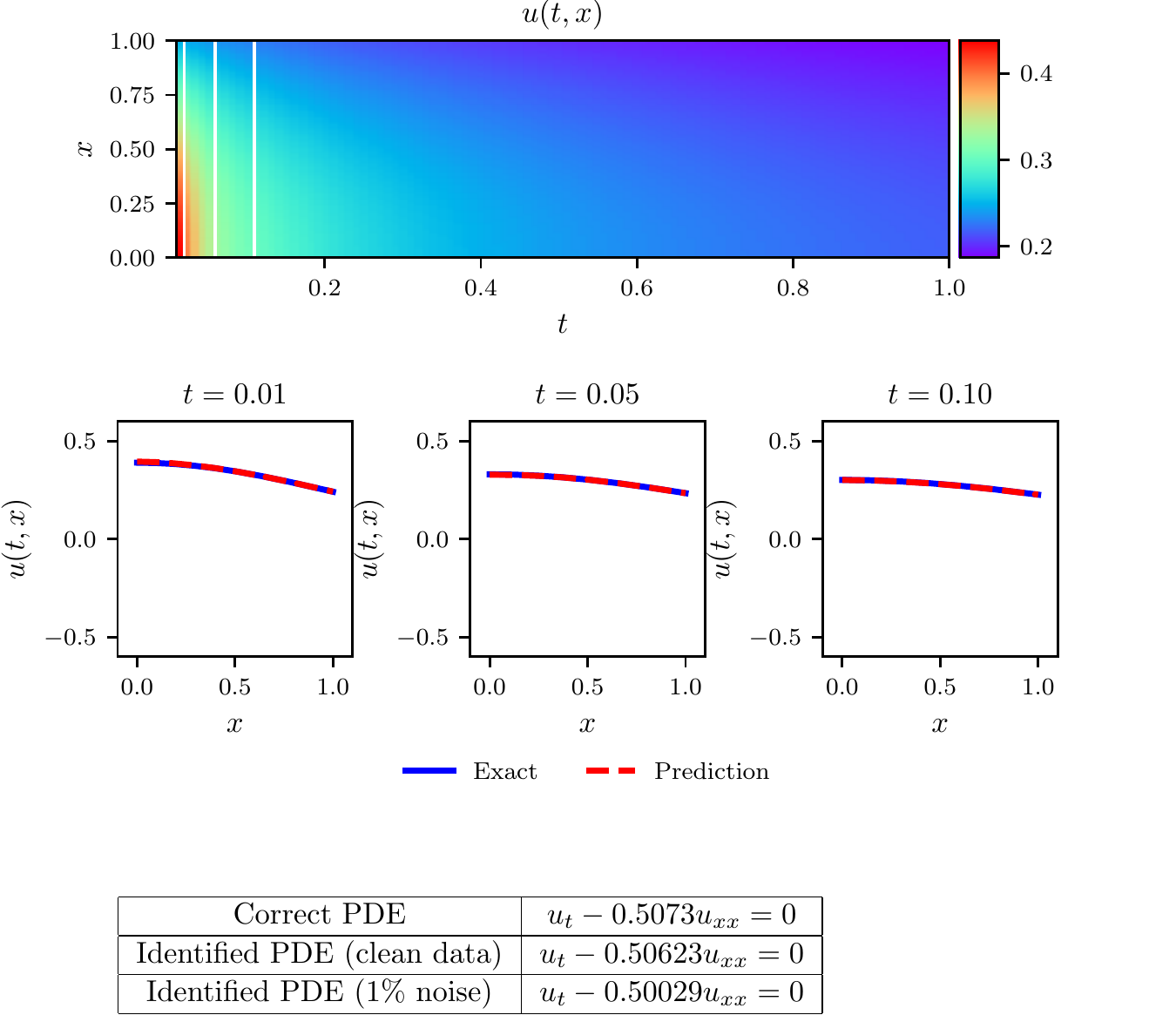}
	\caption{Equation with $\alpha = 0.3$: Top: Solution u(t, x) along with the temporal locations of the three training snapshots. Middle: Training data and exact solution corresponding to the three temporal snapshots depicted by the dashed vertical lines in the top panel. Bottom: Correct partial differential equation along with the identified one obtained by learning $\lambda$.}
	\label{fig21}
\end{figure}
\par In figure \ref{fig21}, we summarize the results for the conformable time-fractional diffusion equation with $\alpha = 0.3$. At the top panel, we exhibit the predicted solution of the equation while its analytical solution can be found in figure \ref{fig6}. At the middle pannel, the graph of the predicted solution and the exact solution at $t=0.01,0.05,0.10$ is drawn. At the bottom pannel, we exhibit the correct PDE and the identified PDE with between clean data and 1$\%$ noise. It is a good news that proposed neural network method can accurately identify the unknown parameter $\lambda$ by making use of our neural network method, even if the training data is corrupted with noise. Concretely, the estimation errors of $\lambda$ is 0.11$\%$ in the case of clean training data. Also, even though the training data is corrupted with 1$\%$ uncorrelated Gaussian noise, the prediction is still robust, returning the error of $\lambda$ is 0.70$\%$.

\noindent{\bf{Example 6.}}
In this example, the equation \eqref{eq1} with $\alpha = 0.8$ is taken into consideration. Here, we set up a training dataset by randomly generating $N = 2000$ points across the entire spatio-temporal domain from the exact solution corresponding to $\lambda = 0.5073$ to illustrate the effectiveness of our method. Then we use $N_{data} = 2000$ training points to train 3021 parameters of the deep neural network that has 9 hidden layers with 20 neurons at each layer by taking advantage of the Adam optimizer and L-BFGS optimizer to minimize the mean square error loss function defined in equation \eqref{eq4} and learn the solution $u(t,x)$ of the conformable time-fractional diffusion equation. At this point, the hyperbolic tangent function is used as our activation function. Afterwards, the network is calibrated to predict the entire solution u(t, x), as well as the unknown parameter $\lambda$. 
\begin{figure}[htb]
	\centering
	\includegraphics[scale=1]{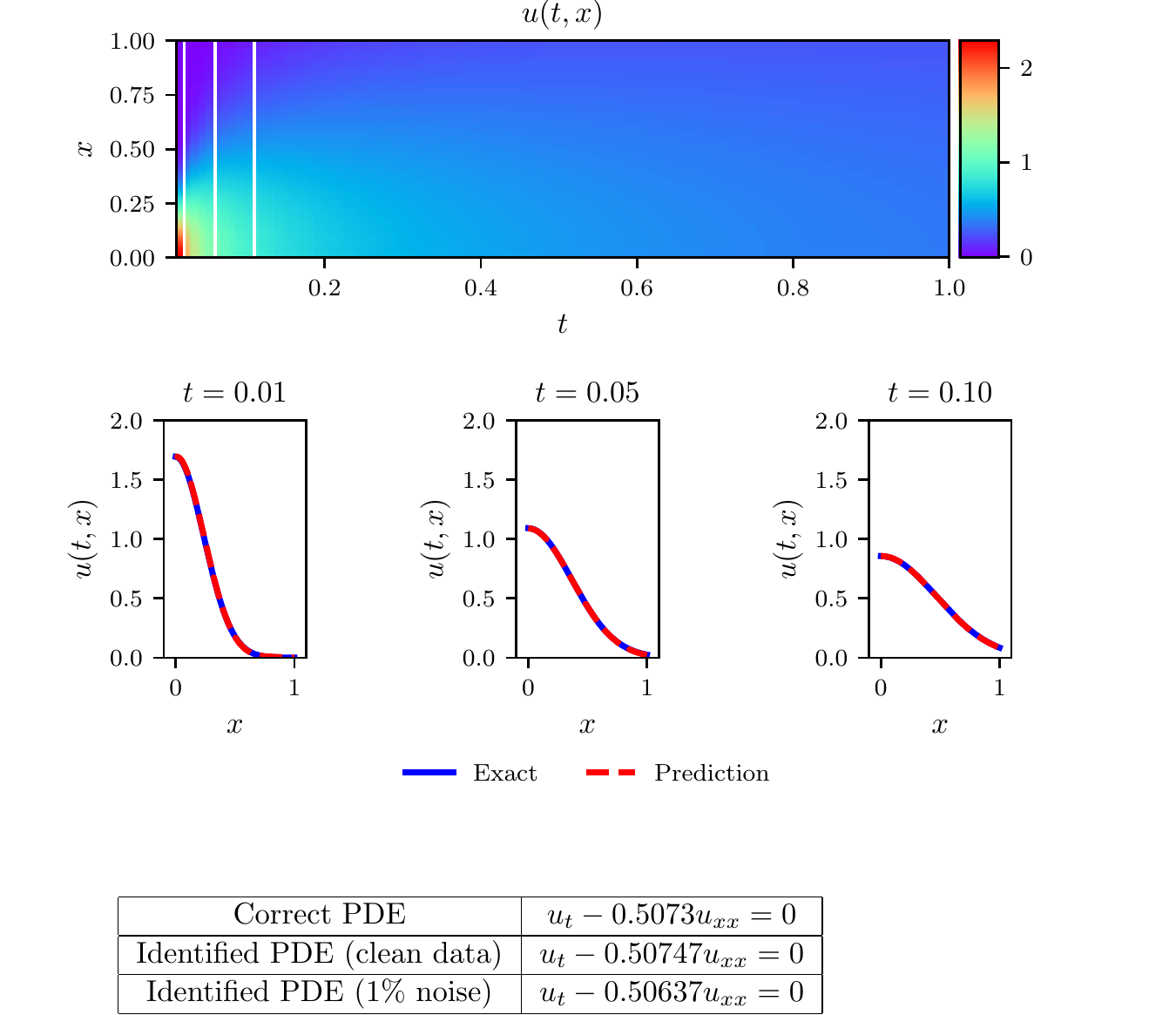}
	\caption{Equation with $\alpha = 0.8$: Top: Solution u(t, x) along with the temporal locations of the three training snapshots. Middle: Training data and exact solution corresponding to the three temporal snapshots depicted by the dashed vertical lines in the top panel. Bottom: Correct partial differential equation along with the identified one obtained by learning $\lambda$.}
	\label{fig22}
\end{figure}
\par In figure \ref{fig22} , we sum up the results for the conformable time-fractional diffusion equation with $\alpha = 0.8$. At the top panel, we show the predicted solution of the equation whose image of analytical solution can be found in figure \ref{fig11}. At the middle pannel, the graph of the predicted solution and the exact solution at $t=0.01,0.05,0.10$ is drawn. At the bottom pannel, we demonstrate the correct PDE and the identified PDE with both clean data and 1$\%$ noise. Through the figure, we observe that our neural network method can accurately identify the unknown parameter $\lambda$, even if the training data is corrupted with noise. Specifically, in the case of clean training data, the estimation errors of $\lambda$ is 0.02$\%$. Despite the training data is corrupted with 1$\%$ uncorrelated Gaussian noise, the prediction is still robust while the error of $\lambda$ is 0.09$\%$.
\par To sum up, our neural network method provides a high accuracy in solving the inverse problem of the conformable time-fractional diffusion equation, and correctly identify unknown parameters regardless of whether the training data is corrupted with noise or not.
It is worth noting that it is incomparable to classical numerical methods.

\section{Conclusion}
In this paper, a neural network method is proposed to study the conformable time-fractional diffusion equation. Since conformable derivative has the nature of chain rule, automatic differential techniques can be applied, and our main idea is to use the PINN method to study the equation. To take this, we have filled in the defect that neural network can not be directly used to solve fractional differential equations in some previous articles. For the forward problem, we train the neural network to obtain predicted solution by using IC/BCs and give three numerical examples to validate the effectiveness. The error of predicted solution and exact solution is discussed. Particularly, when $\alpha$ approaches 1, the simulation effect of the singular solution part of the equation we study is not satisfactory. Therefore, we propose a weighted neural network method that it can constrain the singular solution better, so as to eliminate the influence of the singular solution. Besides, we conduct a systematic study to quantify the effects of different sampling points (training points, collocation points) and neural network structures (layers, the number of neurons in each layer) on the prediction accuracy. The general trend shows that as the total number of training data $N_{u}$ is increased when given a sufficient number of collocation points $N_{f}$ or the number of layers and neurons is increased, the prediction accuracy is increased. For the inverse problem, we use the obtained data to train the neural network by minimizing the loss function we define and give a predicted value of the equation parameter $\lambda$. Three numerical examples are given, we show the correct PDE and the identified PDE with clean data and 1$\%$ noise. We observe that our method can accurately identify the parameters, even if the training data is corrupted with 1\% uncorrelated gaussian noise. 
\par Recently, many scholars have studied partial differential equations by means of neural network, and put forward many good methods. However, for solving the partial differential equations, the rise of deep learning methods also brings a lot of problems. We still have a lot of work to do on the conformable time-fractional diffusion equation. In the future work, we shall mainly consider the following two issues. One one hand, whether different sampling strategies will have some good effects on our neural network method or not is a question of interest to us. It may be a good choice to select sampling points clustered where the solution has large gradient. On the other hand, we are more concerned with how to select the best hyper-parameters of the neural network for the conformable time-fractional diffusion equation. 
\bigbreak

\end{document}